\documentclass{elsarticle}
\usepackage{graphicx}
\usepackage{amssymb}
\usepackage{amsmath}
\usepackage{amsfonts}
\newtheorem{algorithm}{Algorithm}

\newdefinition{definition}{Definition}
\newdefinition{remark}{Remark}

\begin{document}

\begin{frontmatter}

\title{Face-based Selection of Corners in 3D Substructuring}

\author[a]{Jakub~\v{S}\'{\i}stek\corref{cor1}}
\ead{sistek@math.cas.cz}
\author[b]{Marta \v{C}ert\'{\i}kov\'{a}}
\ead{marta.certikova@fs.cvut.cz}
\author[b]{Pavel Burda}
\ead{pavel.burda@fs.cvut.cz}
\author[c]{Jaroslav~Novotn\'{y}}
\ead{novotny@mat.fsv.cvut.cz}
\address[a]{Institute~of~Mathematics, 
Academy~of~Sciences of the~Czech~Republic}
\address[b]{Department of Mathematics,
Faculty of Mechanical Engineering \\ Czech Technical University in Prague}
\address[c]{Department of Mathematics,
Faculty of Civil Engineering \\ Czech Technical University in Prague}
\cortext[cor1]{Corresponding author at: Institute of Mathematics, Academy of Sciences of the Czech Republic, \v{Z}itn\'{a}~25, CZ--115~67 Praha 1, Czech Republic.
Tel.: +420 222 090 710; fax: +420 222 211 638}

\begin{abstract}
In most recent substructuring methods,
a~fundamental role is played by the coarse space. 
For some of these methods (e.g. BDDC and FETI-DP),
its definition relies on a~`minimal' set of coarse nodes (sometimes called corners) 
which assures invertibility of local subdomain problems and also of the global coarse problem. 
This basic set is typically enhanced by enforcing continuity of functions 
at some generalized degrees of freedom, such as average values on edges
or faces of subdomains.
We revisit existing algorithms for selection of corners.
The main contribution of this paper consists of proposing a~new heuristic algorithm for this purpose.
Considering faces as the basic building blocks of the interface, 
inherent parallelism, and better robustness with respect to disconnected subdomains are among features of the new technique.
The advantages of the presented algorithm in comparison to some earlier approaches are demonstrated on three engineering 
problems of structural analysis solved by the BDDC method.
\end{abstract}


\end{frontmatter}

\section{Introduction}
\label{sec:Introduction}

The Balancing Domain Decomposition based on Constraints (BDDC) is a numerically scalable, 
nonoverlapping (substructuring), 
primary domain decomposition method
introduced in 2003 by Dohrmann \cite{Dohrmann-2003-PSC}.
Its algebraic theory developed by Mandel, 
Dohrmann and Tezaur in \cite{Mandel-2005-ATP} demonstrates close relation to FETI-DP introduced 
by Farhat, Lesoinne, and Pierson \cite{Farhat-2000-SDP}: 
the eigenvalues of the preconditioned problem in BDDC and FETI-DP are the same except possibly those equal to 0 and 1 
(see also \cite{Brenner-2007-BFW}, \cite{Li-2006-FBB}, and \cite{Mandel-2007-BFM} for simplified proofs).
These results not only provide the theoretical reasoning for nearly identical performance of BDDC and FETI-DP observed earlier, 
but also imply, that many theoretical results obtained for one method apply readily to the other.

The \emph{coarse space}, defined by constraints on continuity of functions across the interface at \emph{coarse degrees of freedom}, 
is essential for the performance of both methods.
A~historical overview of an evolution of the concept of the coarse space is presented, e.g., by Widlund in \cite{Widlund-2008-DCS} 
and by Mandel and Soused\'{\i}k in \cite{Mandel-2011-CSA}.
The usual basic choice of coarse degrees of freedom is presented by selecting a~set of \emph{coarse nodes} (also called \emph{corners}). 
This set is usually selected to be `minimal' in the sense that it is as small as possible while assuring invertibility of local subdomain problems and of the global coarse problem.
For 2D problems this choice ensures good convergence properties.
However, both methods require additional constraints on some generalized degrees of freedom 
such as average values on edges or faces of subdomains to achieve good efficiency for 3D problems.
This fact was first discovered for FETI-DP: 
experimentally observed in Farhat, Lesoinne, and Pierson \cite{Farhat-2000-SDP} and theoretically 
supported by Klawonn, Widlund and Dryja in \cite{Klawonn-2002-DPF}. 
These observations apply to BDDC through the above-mentioned correspondence between both methods.

A sufficiently robust definition of the coarse space in BDDC and FETI-DP is still not available, 
especially for complex 3D geometries, and existing methods tend to fail for such problems.
Related work on choice of the coarse degrees of freedom has focused on selecting a~small and effective coarse space.
An algorithm for selecting the smallest set of coarse nodes to avoid coarse mechanism is described by Lesoinne in \cite{Lesoinne-2003-FCS}.
Another algorithm, which is already based on pairs of subdomains, was given by Dohrmann in \cite{Dohrmann-2003-PSC}.
This task has been recently further discussed by Bro\v{z} and Kruis in \cite{Broz-2009-HAS} for 2D case.
Klawonn and Widlund in \cite{Klawonn-2003-SCD} and \cite{Klawonn-2006-DPF} minimize a set of more general coarse degrees of freedom 
(like weighted averages over edges and faces) to achieve optimal convergence estimates, introducing the concept of an acceptable path.
Adaptive selection of coarse degrees of freedom based on local estimates using eigenvectors associated with faces is described by Mandel and Soused{\'\i}k in \cite{Mandel-2007-ASF}, and by Mandel, Soused{\'\i}k, and \v{S}\'{\i}stek in \cite{Mandel-2009-ABT}.
In this adaptive approach, which provides additional averages on faces leading to optimal decrease of 
the expected condition number, 
a sufficient number of initial constraints is required between each pair of subdomains as an input.
This assumption is in good agreement with the output of the algorithm proposed in this paper.

While proposing the new algorithm for selection of the basic set of corners is the main contribution of the manuscript,
we further explore the potential of adding more coarse nodes into the coarse problem. 
This approach is technically simple and allows flexible setting of desired approximation. 
It is observed, that by loosening the requirement of `minimal' selection and identifying more interface nodes as corners,
the performance of the BDDC preconditioner may be cheaply but considerably improved.
Numerical experiments on industrial 3D elasticity problems demonstrate the advantages 
of the new corner selecting algorithm in comparison to several earlier approaches.
They also show the fact, that by enhancing the basic set of constraints by additional coarse nodes, 
the computational times might be further reduced.

\section{BDDC method}

In this paper, we study the selection of the initial set of constraints in the context of the BDDC method \cite{Dohrmann-2003-PSC}, 
which is briefly recalled in this section.
However, the main ideas of the paper apply to FETI-DP method as well.

After a discretization of a linearized partial differential equation of elliptic type in a given domain $\Omega$ by means of finite element method (FEM), a system of linear algebraic equations
\begin{equation}
\mathbf{A} \mathbf{x} = \mathbf{f}  \label{eq:Axf}
\end{equation}
with a symmetric positive definite matrix $\mathbf{A}$ and a right-hand side $\mathbf{f}$ is solved for the unknown vector $\mathbf{x}$.
Components of $\mathbf{x}$ represent function values at mesh nodes and they are often called \emph{degrees of freedom}. 
In 3D linear elasticity, there are 3 unknown values of displacement (3 degrees of freedom) at every mesh node.

The first step in the BDDC method is the reduction of the problem to the interface. 
This is quite standard and described in the literature, 
e.g., Toselli and Widlund \cite{Toselli-2005-DDM}: 
the underlying discretized domain $\Omega$ is split into $N$ nonoverlapping subdomains (also called \emph{substructures}) 
$\Omega_i,\ i = 1,\dots,N$ with common interface $\Gamma$, and problem (\ref{eq:Axf}) is reduced to the Schur complement problem with respect to interface
\begin{equation}
\mathbf{S} \mathbf{u} = \mathbf{g}  \label{eq:Sug}
\end{equation}
with a symmetric positive definite matrix $\mathbf{S}$. 
The vector $\mathbf{u}$ now represents the subset of degrees of freedom in $\mathbf{x}$ that correspond to the interface $\Gamma$.
Solution $\mathbf{u}$ of the problem (\ref{eq:Sug}) can be also represented as the minimum of the functional
\begin{equation}
 \frac{1}{2} {\mathbf{u}} ^{\rm T} {\mathbf{S}} {\mathbf{u}} - {\mathbf{u}} ^{\rm T} \mathbf{g} \, \rightarrow \, min  \, , \quad \mathbf{u} \in \widehat{W}  \label{eq:S-min}
\end{equation}
on the space $\widehat{W}$ of unknowns on the interface $\Gamma$. 
The space $\widehat{W}$ can be identified with the space of \emph{discrete harmonic functions}, 
that are fully determined by their values of unknowns on the interface $\Gamma$ and have minimal energy on every subdomain.

The problem (\ref{eq:Sug}) is then solved by the preconditioned conjugate gradient (PCG) method,
for which BDDC acts as the preconditioner. 
The main idea of the BDDC method is shortly described bellow. 
More details, together with connection to FETI-DP, 
can be found in Mandel, Dohrmann and Tezaur \cite{Mandel-2005-ATP} or Mandel and Soused{\'\i}k \cite{Mandel-2007-BFM}.

A preconditioner ${\mathbf{M}}$ for the system (\ref{eq:Sug}) should realize an approximation of ${\mathbf{S}^{-1}}$ such that obtaining a preconditioned residual
$\mathbf{p} = \mathbf{M} \mathbf{r}$ can be considerably easier than solving the original problem (\ref{eq:Sug}). Construction of the BDDC preconditioner is based on the idea that instead of minimising (\ref{eq:S-min}) on the space $\widehat{W}$, which represents solving the system (\ref{eq:Sug}), the minimization is performed on some larger space $\widetilde{W}$ such that $\widehat{W} \subset \widetilde{W}$:
\begin{equation}
 \frac{1}{2} {\mathbf{\widetilde{u}}} ^{\rm T} {\mathbf{\widetilde{S}}} {\mathbf{\widetilde{u}}} - {\mathbf{\widetilde{u}}} ^{\rm T} \mathbf{\widetilde{g}} \, \rightarrow \, min  \, , \quad \mathbf{\widetilde{u}} \in \widetilde{W} , \label{eq:Stilde-min}
\end{equation}
where ${\mathbf{\widetilde{S}}}$ is a symmetric positive definite extension of ${\mathbf{S}}$ to $\widetilde{W}$ and ${\mathbf{\widetilde{g}}}$ is an extension of ${\mathbf{g}}$. The space $\widetilde{W}$ has to be chosen so that the symmetric positive definite extension ${\mathbf{\widetilde{S}}}$ on $\widetilde{W}$ exists. 
At the same time, solving problem~(\ref{eq:Stilde-min}) should be considerably easier than solving the original problem (\ref{eq:S-min}), 
while providing good approximation of the solution of (\ref{eq:S-min}). 
The BDDC preconditioner is then defined as
\begin{equation}
\mathbf{M} = \mathbf{E} {\mathbf{\widetilde{S}}}^{-1} \mathbf{E}^{\rm T} \, ,  \label{eq:M}
\end{equation}
where $\mathbf{E}$ represents a projection from $\widetilde{W}$ onto $\widehat{W}$ realized by a kind of averaging.

\section{Coarse degrees of freedom}

In BDDC, 
the space $\widetilde{W}$ is specified by relaxing the requirement of the continuity of discrete harmonic functions across the interface. 
The functions of $\widetilde{W}$ are required to be continuous only at selected \emph{coarse degrees of freedom}. 
In this paper, we focus on
the simplest choice of coarse degrees of freedom, 
which is a function value at a selected node on the interface. 
Such node is then called \emph{coarse node} or \emph{corner}.
More general coarse degrees of freedom are commented at the end of this section and are considered in computations.

\begin{figure}[htbp]
\begin{center}
\includegraphics[width=37mm]{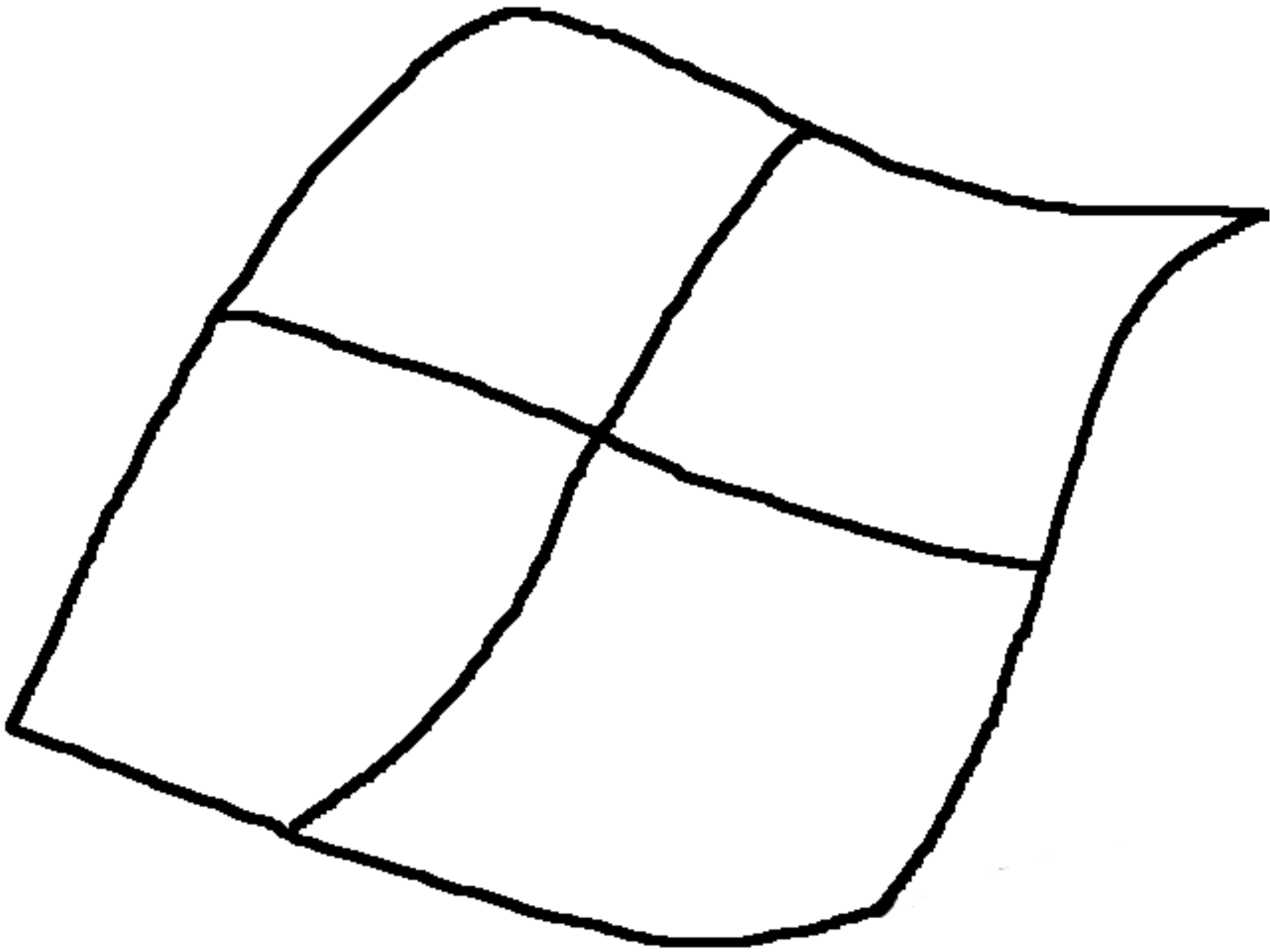} \quad \quad
\includegraphics[width=37mm]{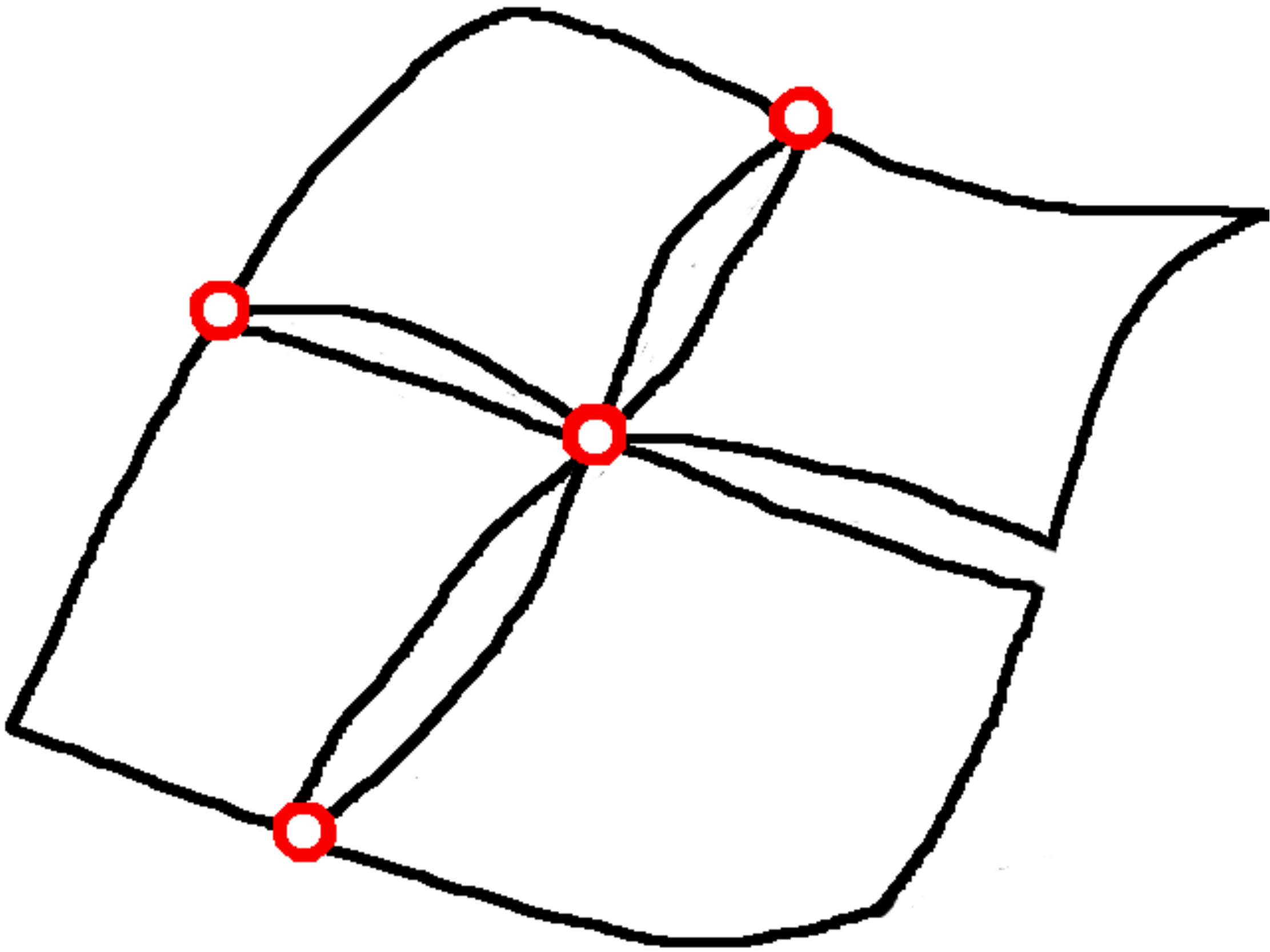}
\includegraphics[width=37mm]{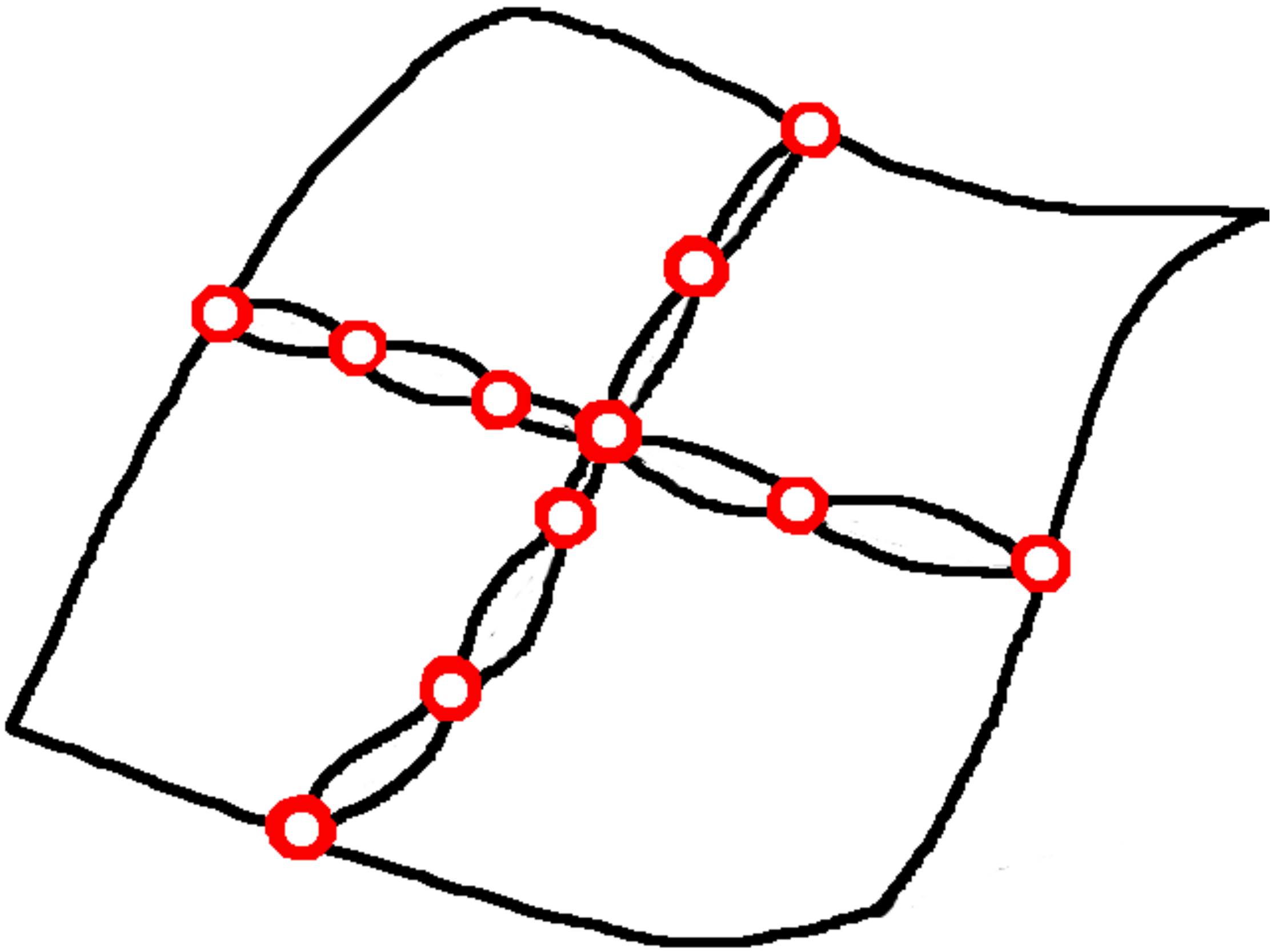}
\end{center}
\hskip 2cm $\widehat{W}$ \hskip 5.5cm $\widetilde{W}$
\caption{A schematic illustration of the continuity constraints: functions from $\widehat{W}$ are continuous across the interface (left), functions from $\widetilde{W}$ are continuous only at corners, marked by circles (centre and right, for two different choices of $\widetilde{W}$).}
\label{fig:Wtilde}
\end{figure}

In terms of mechanics, 
the transition from $\widehat{W}$ to $\widetilde{W}$ can be interpreted as making cuts into the continuous function along the interface, 
leaving the function continuous across the interface only at the corners.
A schematic illustration of the continuity constraints is depicted in Figure \ref{fig:Wtilde}: 
functions from $\widehat{W}$ are continuous across the interface, 
functions from $\widetilde{W}$ are continuous only at selected coarse nodes. 

The space $\widetilde{W}$ can be decomposed as $\mathbf{\widetilde{S}}$-orthogonal direct sum
\begin{equation}
\widetilde{W}=\widetilde{W}_{1}\oplus\dots\oplus\widetilde{W}_{N}
\oplus\widetilde{W}_{C} \, , \label{eq:sum}
\end{equation}
where $\widetilde{W}_{i},\ i = 1,\dots,N,$ are local, subdomain spaces and $\widetilde{W}_{C}$ is the global \emph{coarse space}, 
defined as the $\mathbf{\widetilde{S}}$-orthogonal
complement of all spaces $\widetilde{W}_{i}$, 
i.e. 
\begin{equation}
\mathbf{w}_C \in \widetilde{W}_{C} \  \Leftrightarrow \  \mathbf{w}_C^T \mathbf{\widetilde{S}}\mathbf{w} = 0 \ \ \forall \mathbf{w} \in \widetilde{W}_{i}, \  i = 1,\dots,N.
\end{equation}

Functions from $\widetilde{W}_{i}$ can have nonzero values only in $\Omega_{i}$ except for coarse degrees of freedom. 
They have zero values at coarse degrees of freedom,
and they are fully determined by degrees of freedom on $\Gamma$ and the discrete harmonic extension to interiors of subdomains.
Similarly, functions from $\widetilde{W}_{C}$ are fully determined
by their values at coarse degrees of freedom (where they are continuous) and by the discrete harmonic extension to interiors of subdomains and on the rest of the interface (i.e. everywhere apart from the coarse nodes). Functions from the spaces $\widetilde{W}_{C}$ and $\widetilde{W}_{i}$ are generally discontinuous across $\Gamma$ outside corners.

According to decomposition (\ref{eq:sum}), 
solution of the problem~(\ref{eq:Stilde-min}) can now be split into solution of $N$ local subdomain problems on the spaces $\widetilde{W}_{i}$ and one global coarse problem on the coarse space $\widetilde{W}_{C}$. All these problems are mutually independent and so can be naturally parallelized.

Coarse degrees of freedom have to be selected so that stable invertibility of both the coarse problem and the local problems is assured. 
Important role of the coarse space is to assure scalability by global error propagation over the whole domain.
It was shown that while for 2D elasticity problems the BDDC (or FETI-DP) preconditioner is scalable for coarse space defined by coarse nodes (corners) only, in 3D elasticity problems more general coarse degrees of freedom, such as (weighted) average values over edges and faces, need to be used in order to achieve scalability, see e.g. Toselli and Widlund \cite{Toselli-2005-DDM}.  

Choice of the coarse degrees of freedom has a great impact on the performance of the preconditioner $\mathbf{M}$.
The more coarse degrees of freedom are chosen, 
the more difficult it is to obtain the solution of (\ref{eq:Stilde-min}), 
which, on the other hand, is then closer to the solution of the original problem (\ref{eq:S-min}). 
In the extreme case of selecting all interface nodes as coarse, $\widetilde{W}_{C} \equiv \widetilde{W} \equiv \widehat{W}$, 
coarse problem becomes the original problem (\ref{eq:S-min}) and $\mathbf{M} \equiv \mathbf{S}^{-1}$. 
In the opposite extreme, 
if no coarse degrees of freedom are selected, 
$\widetilde{W}_{C}$ is empty and solution of (\ref{eq:Stilde-min}) splits to $N$ local problems only, 
some of which might not be invertible. 
Thus, the optimal choice of the coarse space lies somewhere in-between.

\section{Geometry and selection of the coarse space in 3D}

The interface $\Gamma$ in 3 dimensions can be specified as a set of nodes belonging to at least two subdomains (subdomains are considered as closed sets). 
It consists of subdomain faces, edges and vertices. 
While there is an intuitive geometric notion what these three entities mean in a~simple case of a cubic domain divided into cubic subdomains, 
there is no unique exact classification in more general case of domain with complicated geometry and subdomains obtained, e.g., 
by a graph partitioning tool.
We adopt the classification presented by Klawonn and Rheinbach in \cite{Klawonn-2007-RFM} and use it in a slightly simplified form, 
which does not assume knowledge of boundary of the domain and is easy to implement:

\begin{definition}
\label{def:classification}
\quad 
\begin{itemize}
	\item a  \emph{face}   contains all nodes shared by the same two subdomains,
   	\item an \emph{edge}   contains nodes shared by the same more than two subdomains,
   	\item a  \emph{vertex} is a degenerated edge with only one node.
\end{itemize}
\end{definition}

Then every node of the interface belongs to just one of the entities defined above. 
Two subdomains are called \emph{adjacent} if they share a face.

\begin{figure}[htbp]
\hbox{
\centering
\vbox{
	\hsize=40mm
	\centering
	\includegraphics[width=20mm]{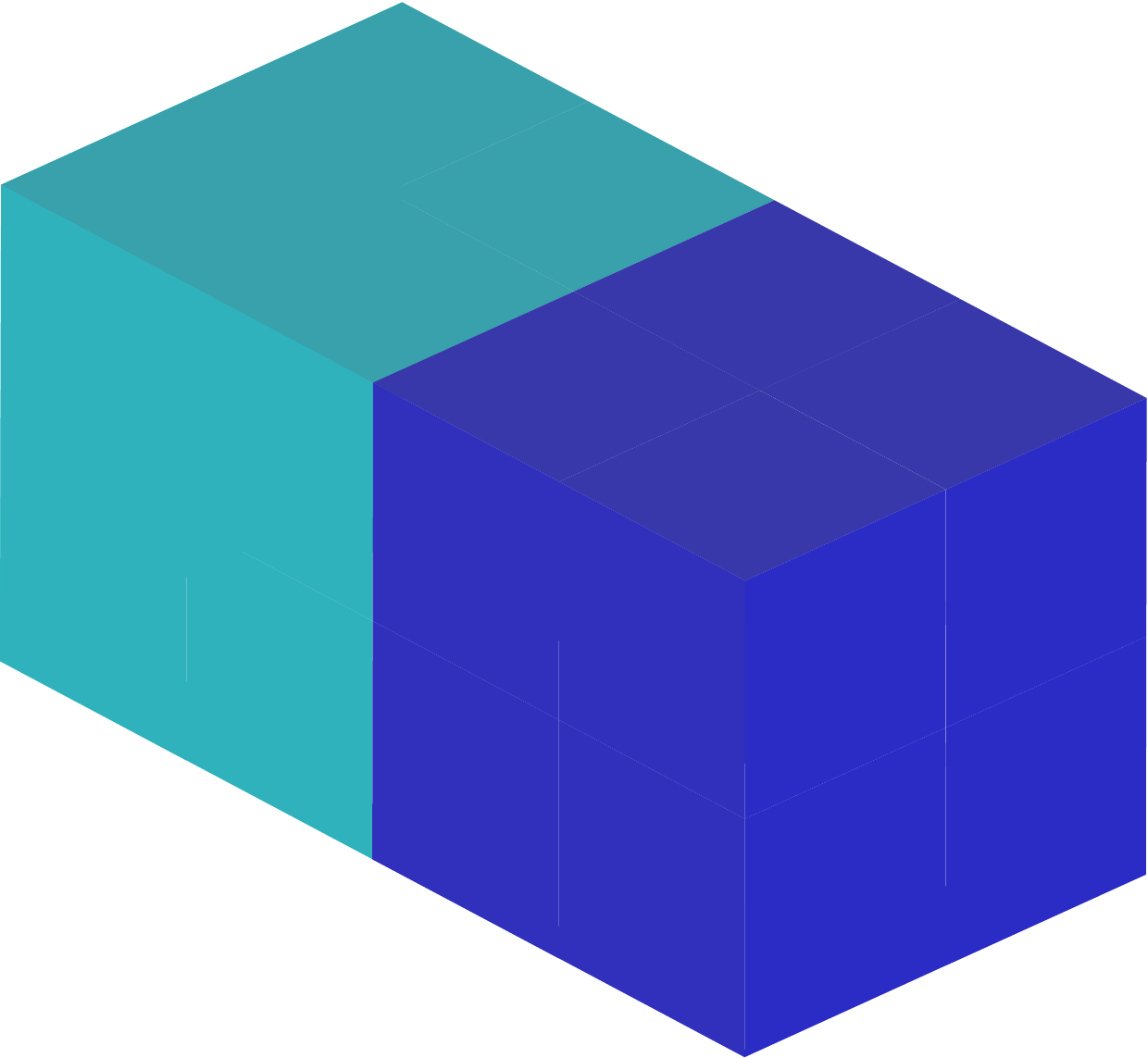}
}
\vbox{
	\hsize=40mm
	\centering
	\includegraphics[width=25mm]{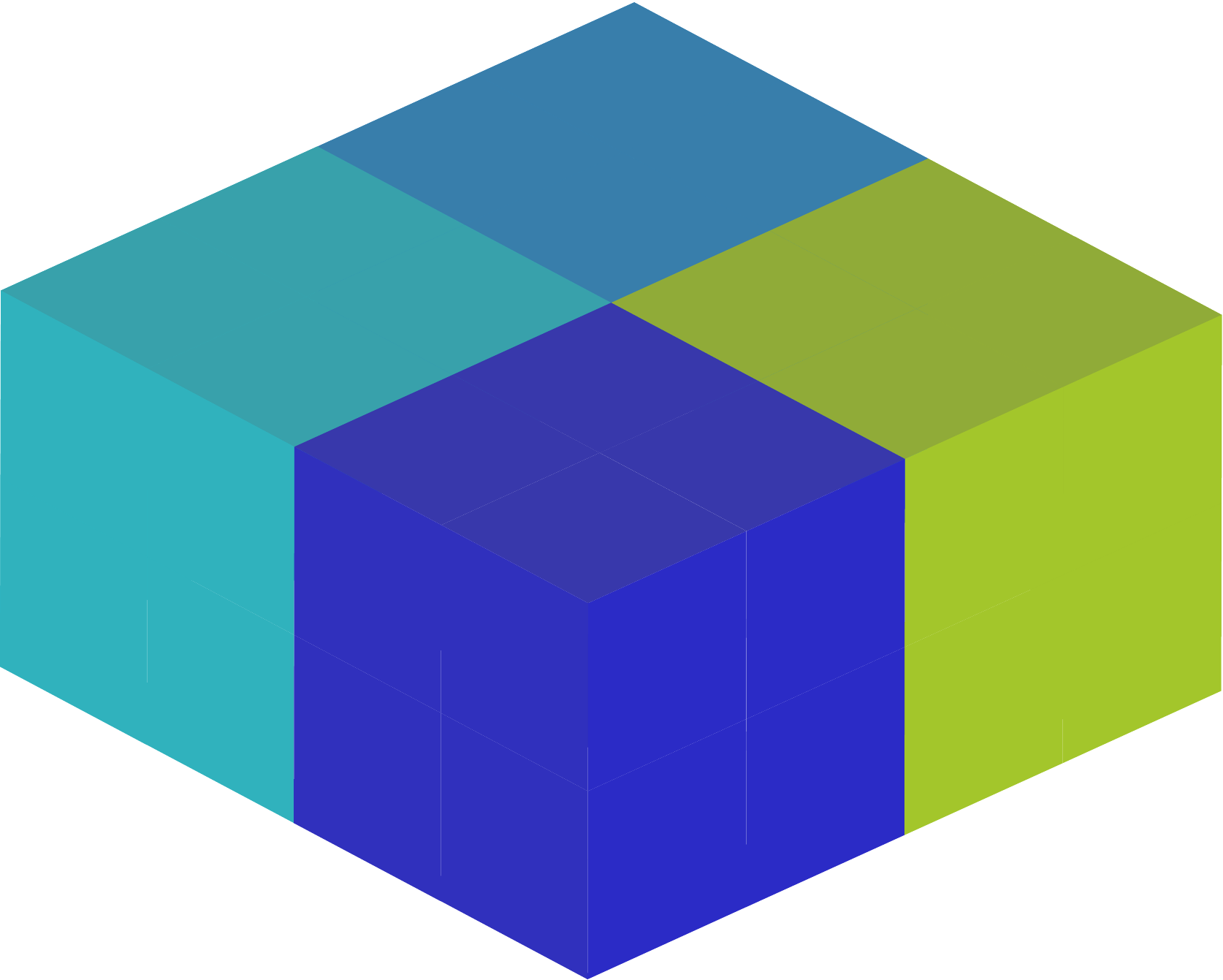}
}
\vbox{
	\hsize=40mm
	\centering
	\includegraphics[width=17mm]{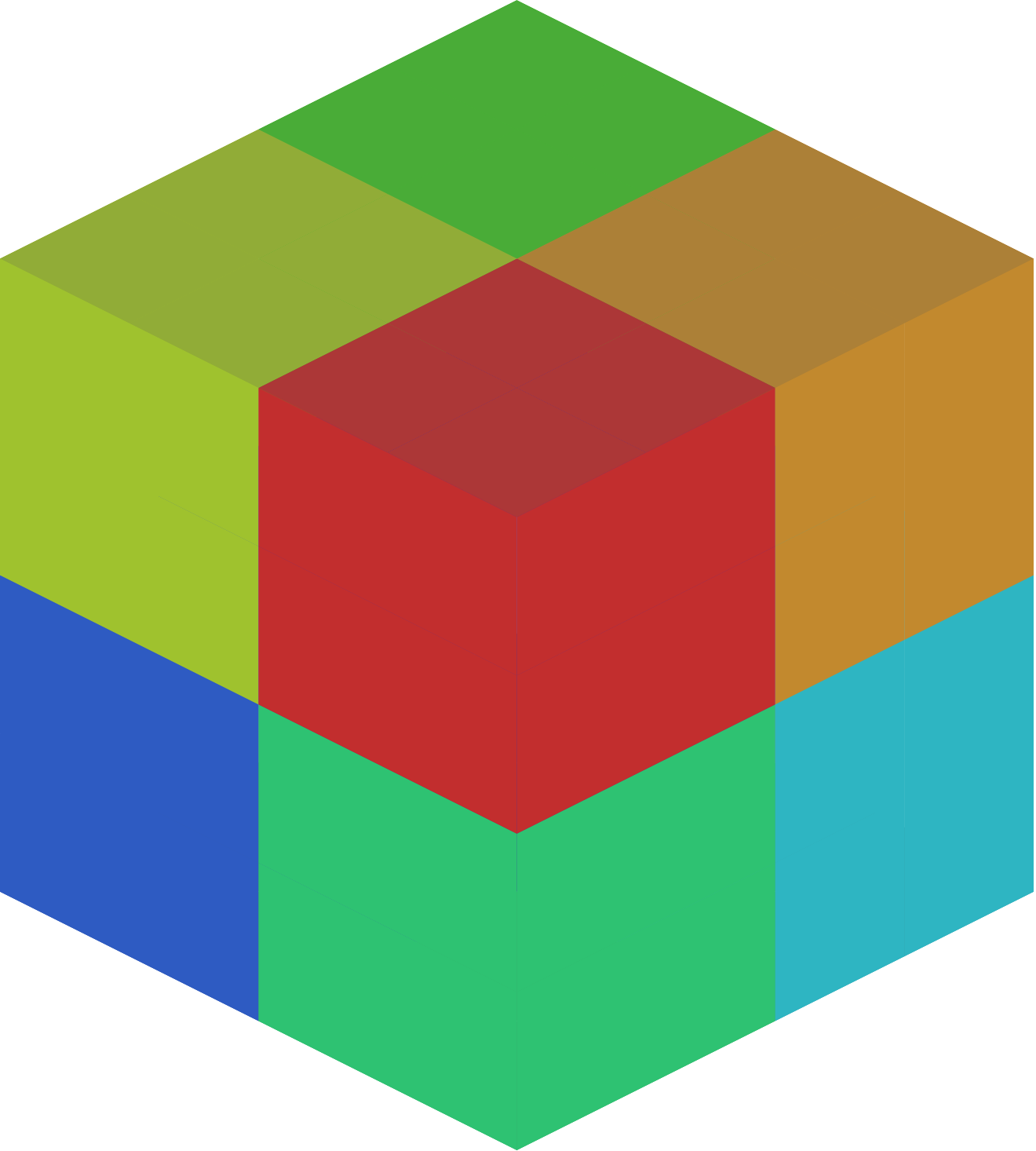}
}
}
\vskip 3mm
\hbox{
\centering
\vbox{
	\hsize=40mm
	\centering
	1 face
	\vspace{5mm}
}
\vbox{
	\hsize=40mm
	\centering
	4 faces \\
	1 edge
	\vspace{3mm}
}
\vbox{
	\hsize=40mm
	\centering
	12 faces \\
	6 edges \\
	1 vertex
}
}
\caption{Examples of classification of the interface nodes as faces, edges and vertices according to \emph{Definition 1}.}
\label{fig:cubes}
\end{figure}

However, this classification does not reproduce our intuition in the case of cubic subdomains, 
as can be seen in Figure \ref{fig:cubes}: 
for instance the interface of a~domain consisting of two cubic subdomains has neither vertex nor edge, 
just one face (the left case in Figure \ref{fig:cubes}). 
Different definitions of faces and edges are discussed by Klawonn and Rheinbach in \cite[Section 2]{Klawonn-2006-PID}.

In practice, 
there are often not enough vertices, edges, or faces for satisfactory number of constraints. 
We have found it useful to introduce one additional entity:
\begin{definition}
\quad 
\begin{itemize}
 \item a  \emph{corner} is any interface node selected as coarse.
\end{itemize}
\end{definition}

In implementations of the BDDC method, 
it is often customary to distinguish between the following two kinds of constraints on continuity across interface.
\vskip 1mm
\textbf{Node constraints - corners}\\
The most obvious choice of coarse degrees of freedom are node constraints (at corners). 
The basic choice is a set of corners, 
that assures invertibility of local subdomain problems and also the global coarse problem.
This is often put as a~requirement on their selection (e.g. in \cite{Dohrmann-2003-PSC},  \cite{Sistek-2010-BFS}).

Although vertices provide a good initial set of corners, 
they often do not suffice for assuring invertibility of subdomain problems and/or of the coarse problem
(cf. Figure \ref{fig:cubes}), and other constraints need to be added.
When other nodes are selected as corners, 
they have to be excluded from corresponding faces or edges, 
so that every interface node is either a corner, or belongs to a face or an edge.

\begin{figure}[htbp]
\begin{center}
\includegraphics[width=39mm]{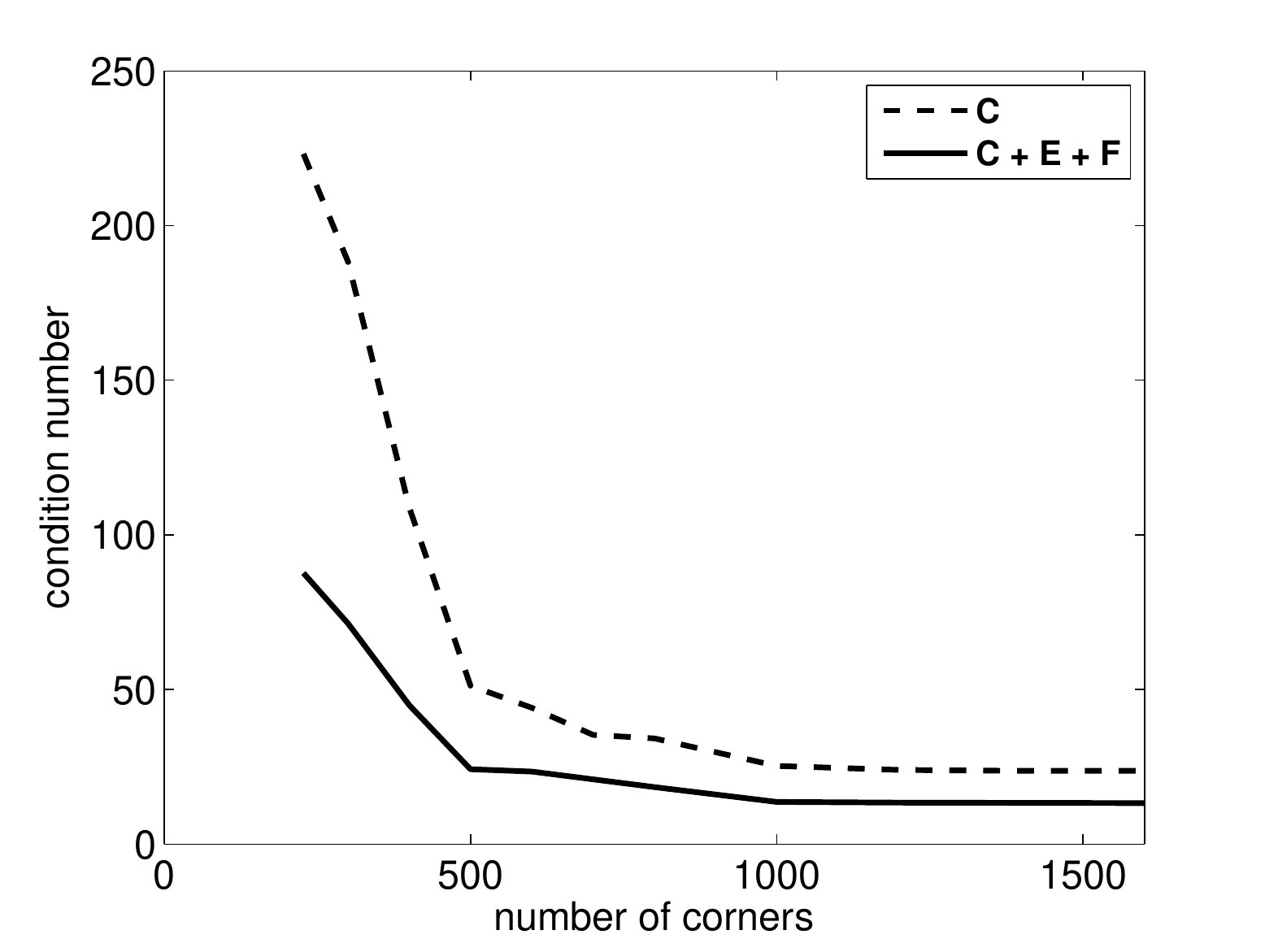}
\includegraphics[width=39mm]{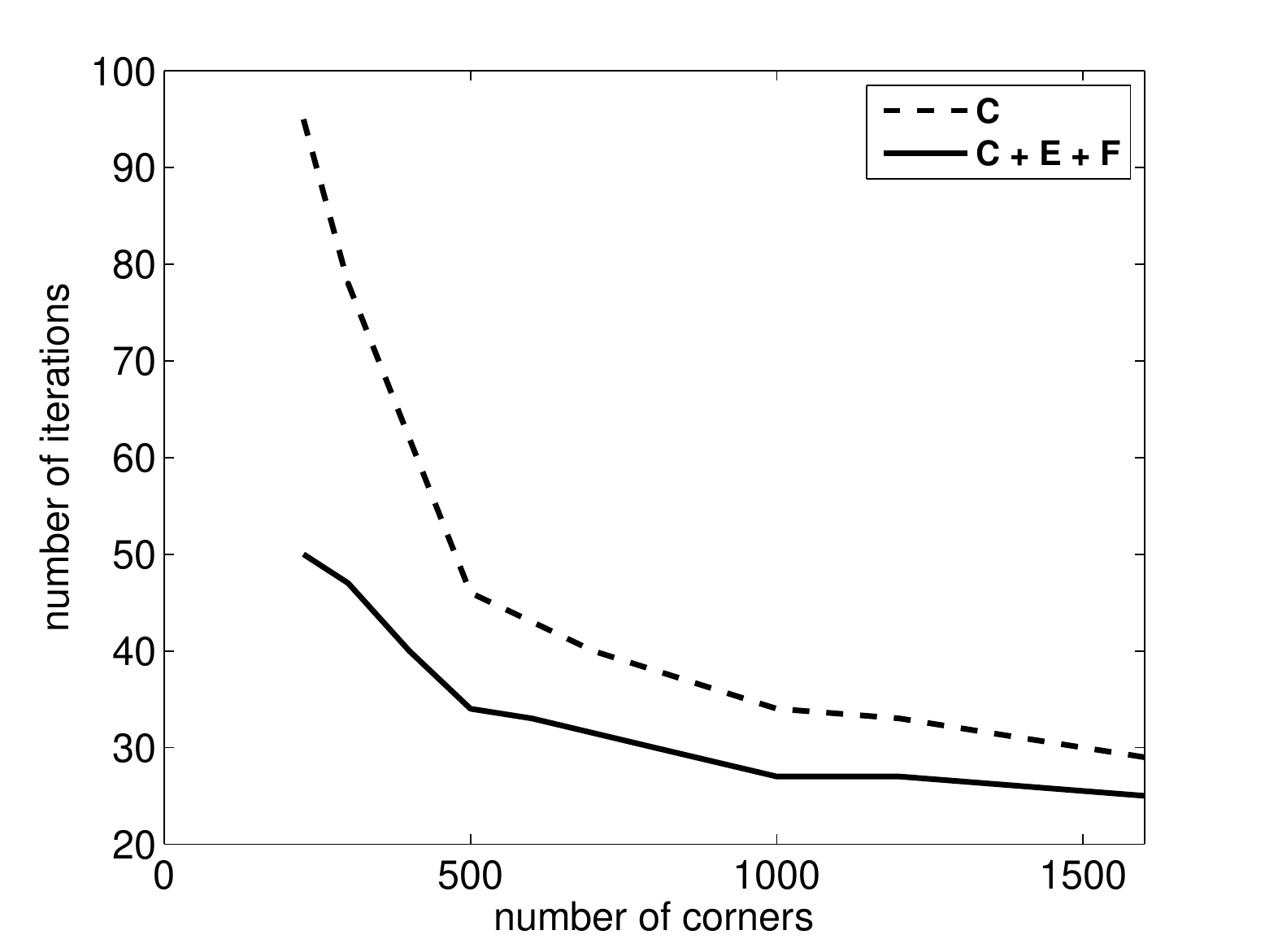}
\includegraphics[width=39mm]{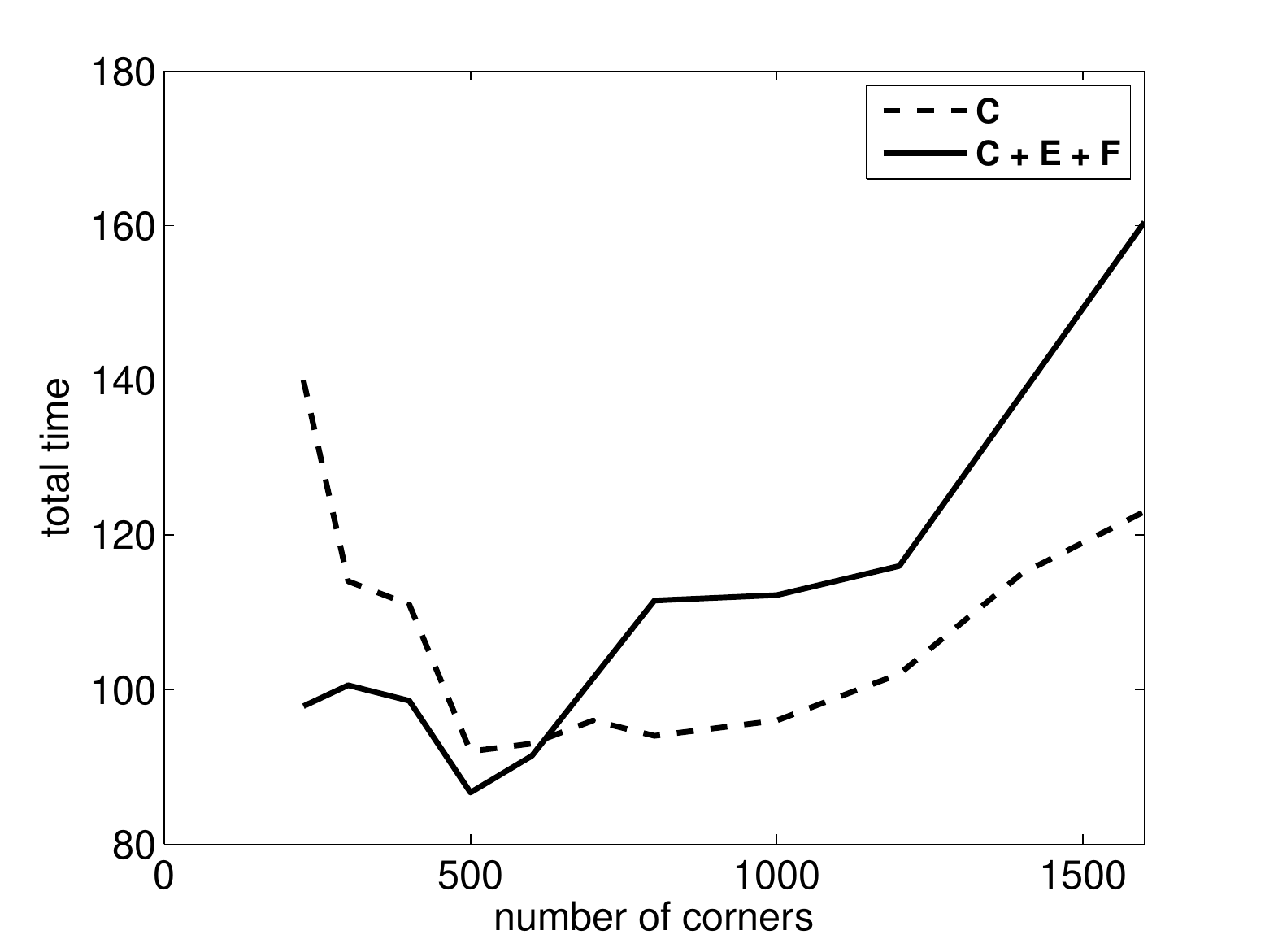}
\end{center}
\caption{Typical dependence of the condition number (left), the number of iterations of the PCG (centre), and the total computational time (right) on the number of corner constraints. Dashed line - corner constraints only, full line - corner constraints and all face and edge averages. Hip joint replacement, 33 186 quadratic elements, 36 subdomains.}
\label{fig:celekt_condition}
\end{figure}

Corner constraints are not as efficient as constraints on averages on edges or faces, 
nevertheless they can be used as a substitute for these constraints, 
if enough corners are employed. 
Figure \ref{fig:celekt_condition} left illustrates the typical dependence of the condition number of the preconditioned problem on number of corners randomly selected from the interface, 
starting from some basic set. 
For small numbers of corners, 
we can observe poor performance of the preconditioner even though all system matrices are invertible. 
Then, after a typical sudden drop, the condition number improves only slightly with adding more corners.
Number of iterations reproduces this dependence, see Figure \ref{fig:celekt_condition} (centre).

Improving convergence by adding more corners leads to a larger coarse problem than adding averages on faces or edges. 
On the other hand, its implementation is straightforward and its scaling is easy to maintain.

For 2D problems, the basic set of corner constraints already ensures good convergence properties. 
Although an efficient BDDC method for 3D elliptic problems requires also constraints on some generalized degrees of freedom, 
such as average values on edges or faces of subdomains described below, 
for many industrial problems this simple approach also leads to satisfactory results. 
\vskip 1mm
\textbf{Constraints on averages over edges and faces}\\
General coarse degrees of freedom can be constructed as any linear combinations of function values at nodes belonging to one face or one edge. 
This type of constraints is required for both BDDC and FETI-DP methods in three dimensions,
if one expects the optimal polylogarithmic bound on condition number $\kappa$ of the preconditioned operator
\begin{equation}
\kappa(\mathbf{M}\mathbf{S}) \leq const.\left(  1+\log\frac{H}{h}\right)  ^{2},
\label{eq:polylog-bound}
\end{equation}
where $H$ is the subdomain size and $h$ is the finite element size (see \cite{Klawonn-2002-DPF}).

One of the standard choices is an arithmetic average over unknowns separately for each component of displacement leading to three constraints for 3D elasticity. 
We have tested this standard choice applied to all edges, to all faces, or both.
More sophisticated methods of weighted averaging were developed, e.g., by Klawonn and Widlund \cite{Klawonn-2006-DPF},
by Mandel and Soused{\'\i}k \cite{Mandel-2007-ASF}, or recently by Mandel, Soused{\'\i}k, and \v{S}\'{\i}stek \cite{Mandel-2009-ABT}.

\section{Selection of the basic set of corners}
\label{sec:selection}

In this section, 
we concentrate on the selection of the basic set of corners that leads to positive definiteness of matrix 
$\mathbf{\widetilde{S}}$ in (\ref{eq:Stilde-min}).
This task is equivalent to assuring invertibility of both local subdomain problems and the global coarse problem only by corner constraints,
which is often required by implementations (cf. \cite{Dohrmann-2003-PSC}, \cite{Sistek-2010-BFS}).
Therefore, we investigate selection of corners independently of enforcing constraints on general coarse degrees of freedom.

From the mechanical point of view, the question of assuring invertibility of local subdomain problems corresponds to enforcing
enough boundary conditions on a~body to fix rigid body modes, 
with subdomain playing the role of the body. 
This goal is easily attained by selecting three nodes (not in a line) of the interface of a subdomain as corners.

\begin{figure}[htbp]
\begin{center}
\includegraphics[width=60mm]{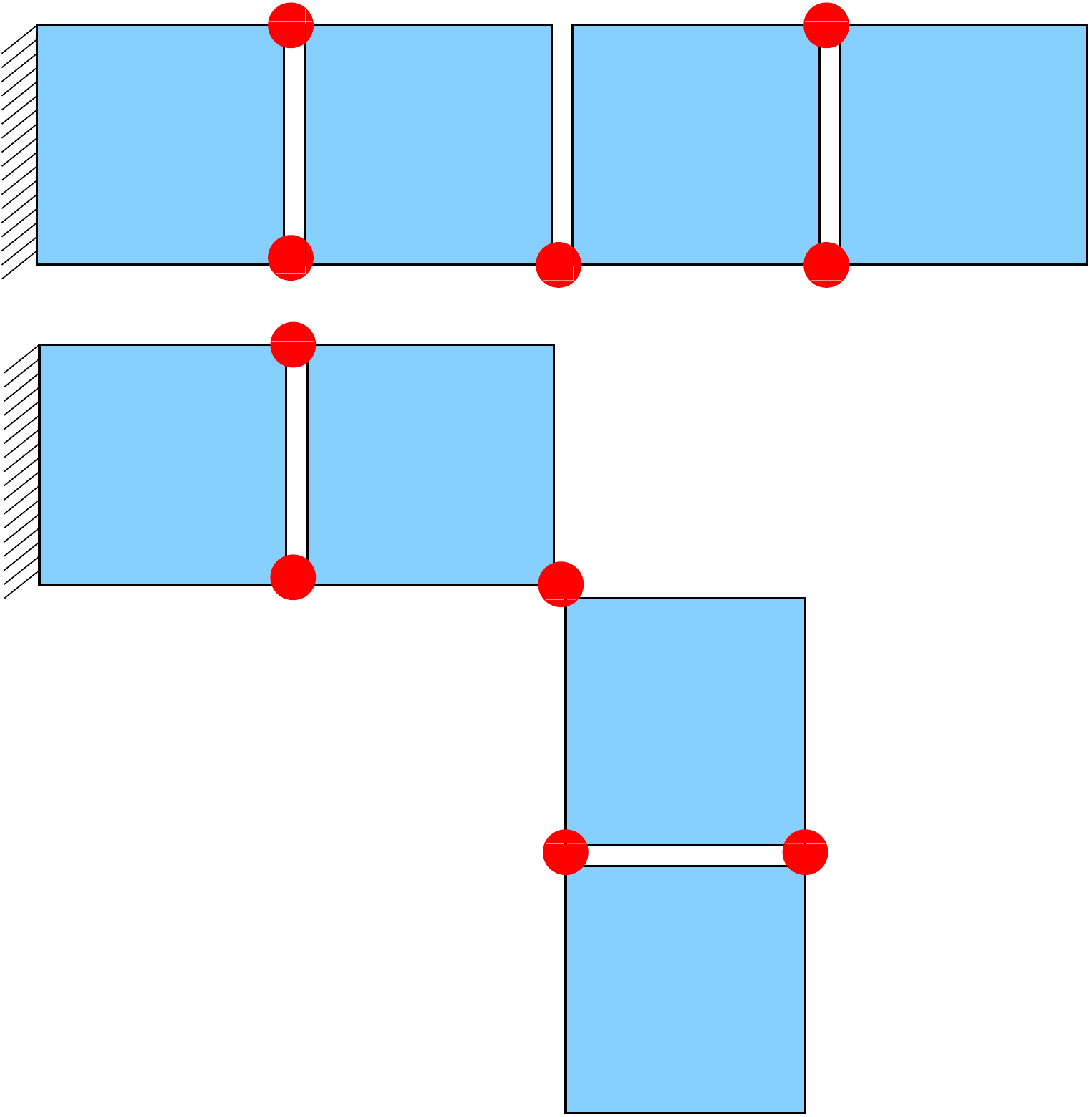}
\caption{A 2D example of mechanism in the coarse problem for serial division into four subdomains, 
red dots denote corners.}
\label{fig:mechanism}
\end{center}
\end{figure}

It turns out, that assuring invertibility of the coarse matrix is the more difficult task, 
since selection with respect to subdomain problems only 
may still lead to \emph{mechanisms} in the coarse problem (see \cite{Lesoinne-2003-FCS}).
To see this, one can simply think of a domain divided into subdomains in a linear fashion.
Figure \ref{fig:mechanism} illustrates this on a 2D case, where two corners for each subdomain are sufficient for invertibility of subdomain stiffness matrices.

An algorithm attempting to select the smallest set of coarse nodes to avoid coarse mechanisms was given by Lesoinne in \cite{Lesoinne-2003-FCS}.
Minimization of the number of corners is obtained mainly by favouring already selected corners.
Thus, the approach is serial in its nature.

Another algorithm for selecting corners was described already by Dohrmann in~\cite{Dohrmann-2003-PSC}.
It is based on the investigation of all possible neighbourings between substructures and selecting three corners from each such set,
that maximise the area of a triangle with corners at its vertices.
However, this algorithm is based on an incomplete classification of interface into vertices, edges, and faces,
and it does not distinguish between the last two groups.
Also this algorithm favours already selected corners by selecting vertices on the interface as the initial vertices of the triangle to be maximised.
Nevertheless, it has provided a good starting point for the new algorithm proposed here.

The third algorithm, which is based on selection of corners along edges, was described in \cite{Sistek-2008-FEM}.
This idea is inspired by the definition of corners as end-points of edges by Klawonn and Widlund \cite{Klawonn-2003-SCD}.
Although it was successfully used by our group in a number of practical computations, 
it may fail to produce a mechanism-free coarse problem in the case of divisions where no edges are present (cf. the leftmost case in Figure \ref{fig:cubes}).

The aim to select a low number of corners inherent to all these algorithms is motivated by the fact, 
that low number of corners results in a small size of the matrix of the coarse problem and its cheap factorization.
However,
it has been observed on a number of experiments (e.g. \cite{Sistek-2010-BFS}, also Section \ref{sec:results} in this paper)
that this motivation may be misleading,
and in fact, larger sets are preferable for the performance of the preconditioner often resulting in much lower number of PCG iterations.
It has been also shown, 
that using more corners may lead to a considerable reduction of the computational time 
in spite of the longer time spent in factorization of the larger matrix of the coarse problem,
even in the case of considering averages on edges and faces.

Based on these observations and experience with the algorithms, 
we see several ideas that the new proposed algorithm should reflect: 
\begin{enumerate}[(i)]
\item selection with respect to faces (by Definition \ref{def:classification}) as these are the basic building blocks of interface in 3D structures (Figure \ref{fig:cubes}),
\item provide larger set of corners than the previous algorithms as this usually leads to much better preconditioning,
\item independence of selection subdomain by subdomain and of order of going through subdomains (better parallelization).
\end{enumerate}

Points (ii) and (iii) are attained simply by not favouring already selected corners
and selecting optimal distribution of at least three corners between each pair of substructures sharing a face, i.e. adjacent substructures, independently.

Let us now present an algorithm satisfying these requirements.
For this, denote the set of faces of subdomain $\Omega _i$ as $\mathcal{F}({\Omega_i})$ and recall that $N$ denotes the number of subdomains.
A face $\mathcal{F}_{ij}$ between subdomains $\Omega_i$ and $\Omega_j$ is present in both sets $\mathcal{F}({\Omega_i})$ and $\mathcal{F}({\Omega_j})$.

\begin{algorithm} [Selection of corners for 3D elasticity problems]
\label{alg:selection}
\quad
\begin{enumerate}
   \item Classify interface according to Definition \ref{def:classification} and use all vertices as corners.
   \item \textbf{For} subdomain $\Omega_i,\ i = 1,\dots,N$,

      \quad \textbf{For} face $\mathcal{F}_{ij} \in \mathcal{F}({\Omega_i}),\ j = 1,\dots,\mathrm{size}(\mathcal{F}({\Omega_i}))$,
   \begin{itemize}
       \item find the set of \emph{all} nodes shared with the adjacent subdomain (generally larger set than the face under consideration, 
             as it may contain also edges and/or vertices),
       \item construct a graph of nodes of this set with connections defined by elements, and detect components of this graph

       \item \textbf{For} each such component, select (in 3D) three corners by:
       \begin{itemize}
           \item[(a)] pick an arbitrary node of the subset,
           \item[(b)] find the first corner as the most remote node from the arbitrary node,
           \item[(c)] find the second corner as the most remote node from the first corner,
           \item[(d)] find the third corner as the node maximising the area of the triangle,
       \end{itemize}
      \textbf{end},
   \end{itemize}
      \quad \textbf{end},

      \textbf{end}.
   \item Select corners as the union of vertices and face-based selection above.
   \item Remove selected corners from edges and faces.
\end{enumerate}
\end{algorithm} 

The algorithm assures that at least three corners are selected in an optimal way with respect to each face.
This situation is often not obtained by favouring already selected corners, 
since corners optimally distributed for one pair of subdomains may be far from optimal distribution with respect to another pair.
Presented algorithm is also much simpler for parallelization than algorithms favouring already selected corners, 
since communication is needed only at the end of the selection to synchronise locally detected corners. 
It typically provides more corners than algorithms mentioned above,
which we consider as an advantage rather than a drawback.

\begin{remark} 
\label{remark:favouring} 
A modification of Algorithm \ref{alg:selection} favouring already selected corners is simply possible by entering the face-based selection in any point (a), (b), (c), or (d),
depending on how many corners are already selected between adjacent substructures.
This modification leads to selection that is very similar to the algorithm by Dohrmann in \cite{Dohrmann-2003-PSC}.
In our experience, this modification, referred to as `minimal', 
leads to lower number of corners, but also usually to worse results (some of them are presented in Section \ref{sec:results}).
Thus, we recommend using the (`full') version as stated by Algorithm \ref{alg:selection}.
\end{remark} 

\begin{remark} 
\label{remark:2D} 
A modification of Algorithm~\ref{alg:selection} for 2D problems (where no edges are present) or topologically 2D problems 
(such as for shell elements in 3D) is simply possible by finishing the face-based selection with point (c).
\end{remark} 

\begin{figure}[htbp]
\begin{center}
\includegraphics[width=55mm]{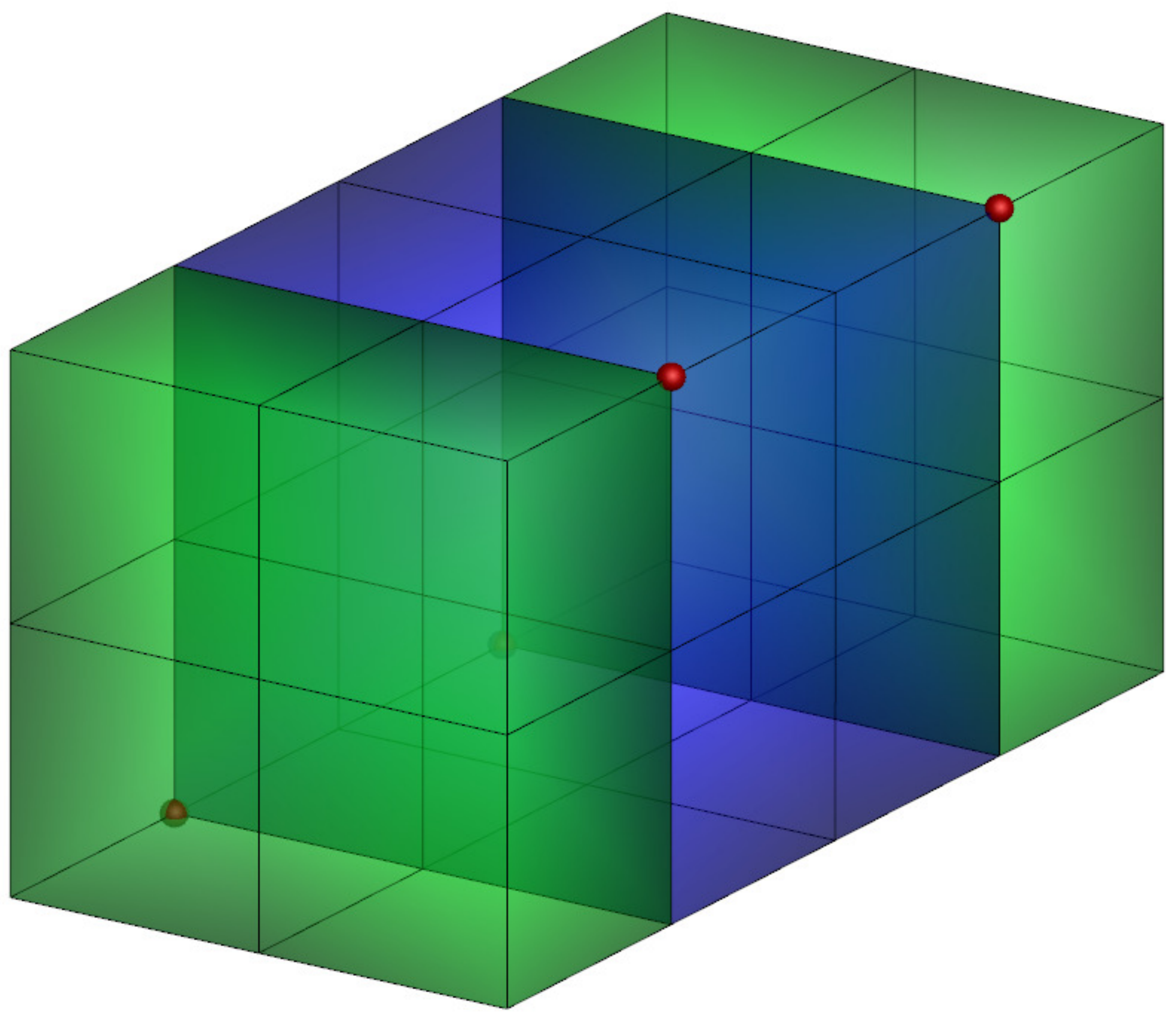}
\includegraphics[width=55mm]{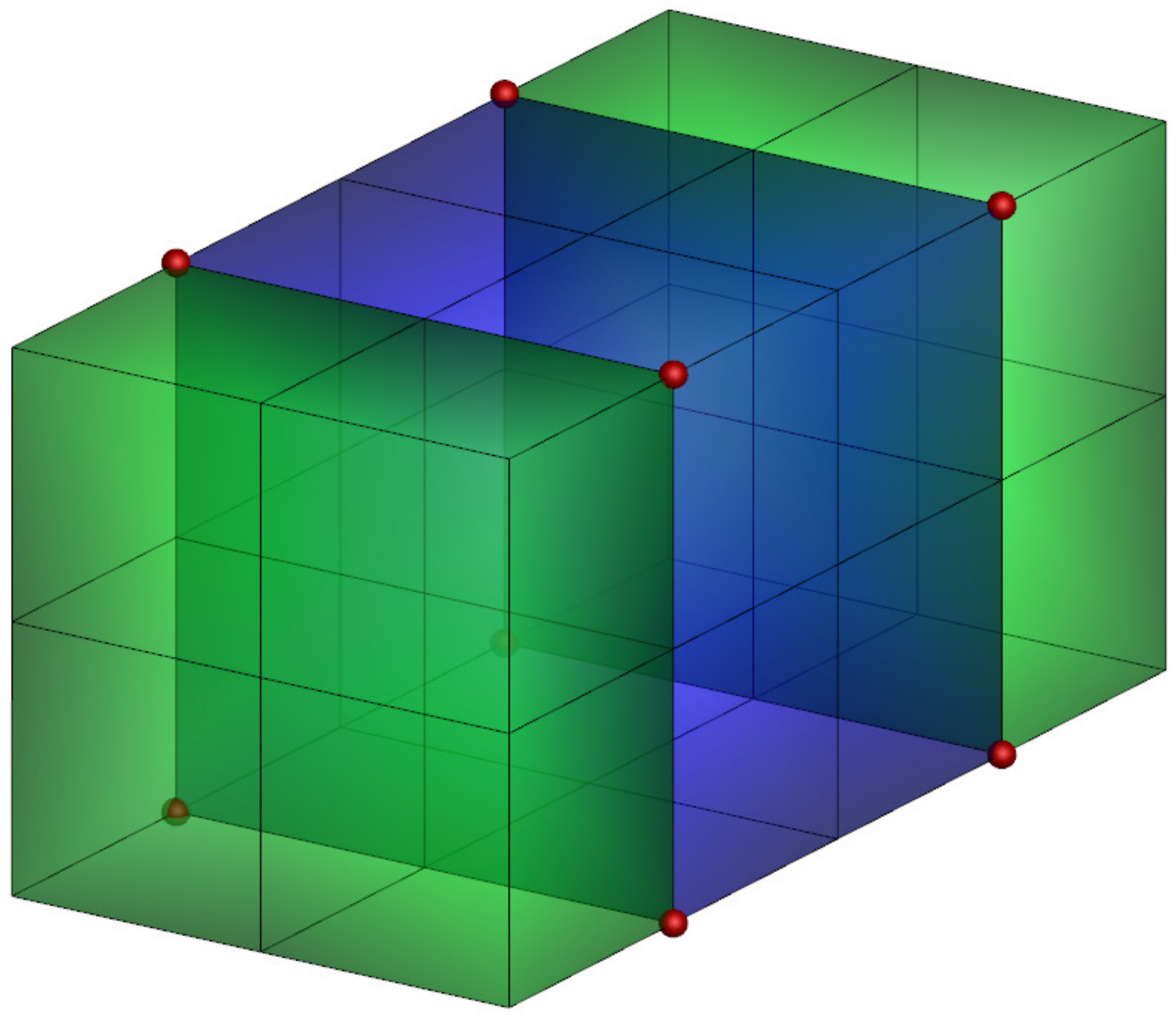}
\end{center}
\caption{Example of a problem with one subdomain disconnected. Four corners obtained by algorithm without detection of components (left), and eight corners obtained with detection of components of interface (right).}
\label{fig:disconnected}
\end{figure}

\begin{remark} 
\label{remark:components} 
Detection of components is aimed at problems divided into subdomains by graph partitioners, such as METIS. 
These programs typically provide divisions well balanced with respect to size of subdomains, 
but often with some subdomains disconnected. 
Such divisions present a~challenge for many existing domain decomposition methods.
With the detection of components, the algorithm is able to detect many of such discontinuities, and fix each component independently. 
The BDDC method is then able to proceed with computation, keeping such subdomains disconnected, 
thus preserving the suggested balance of load. 

We show the power of this detection on a simple problem of an elastic beam consisting of two subdomains, 
one of which is wedged in the other as in Figure \ref{fig:disconnected}. 
On the left-hand side of the figure, corners selected without the detection of components are presented. 
In this case, both cuts are handled as a~single interface, and the search of triangle with maximal area does not succeed.
Resulting configuration has a~mechanism in the coarse problem and, consequently, BDDC method fails.

On the right-hand side of the figure, corners obtained with the component detection enabled are shown. 
Now the optimal triangle is sought at each of the cuts, which leads to a~mechanism-free configuration of corners, 
and the BDDC method converges in four iterations.
\end{remark} 

\section{Implementation}

The BDDC method has been implemented on top of common components of
existing finite element codes, 
namely the frontal solver and the element stiffness matrix generation. 
Such implementation requires only a minimal amount of additional code. 
In our case,
most of the program is written in Fortran 77, with some parts in Fortran 90. 
The MPI library is used for parallelization.

The implementation relies on the separation of node constraints and enforcing the rest by Lagrange multipliers, 
as suggested already in Dohrmann~\cite{Dohrmann-2003-PSC}. 
One new aspect of the implementation is the use of reactions, which come naturally from the frontal solver, to avoid custom coding.
An~external parallel multifrontal solver MUMPS~\cite{Amestoy-2000-MPD} is used for the solution of the coarse problem, 
instead of the serial frontal solver, 
as dimension of the coarse space could become a bottleneck.

Detailed description of the implementation can be found in \cite{Sistek-2010-BFS},
and some more experiments were presented in \cite{Sistek-2008-FEM}.

Recently, the proposed selection of corners has been implemented into the parallel solver,
and the natural parallelism of the algorithm is fully exploited.

\section{Numerical results}
\label{sec:results}

Presented numerical results were computed on SGI Altix~4700 computer with 1.5~GHz Intel Itanium~2 processors (OS Linux) in Czech Technical University Supercomputing Centre, Prague.
For decompositions, we use the METIS graph partitioner \cite{Karypis-1998-FHQ}.

\begin{figure}[htbp]
\begin{center}
\includegraphics[width=55mm]{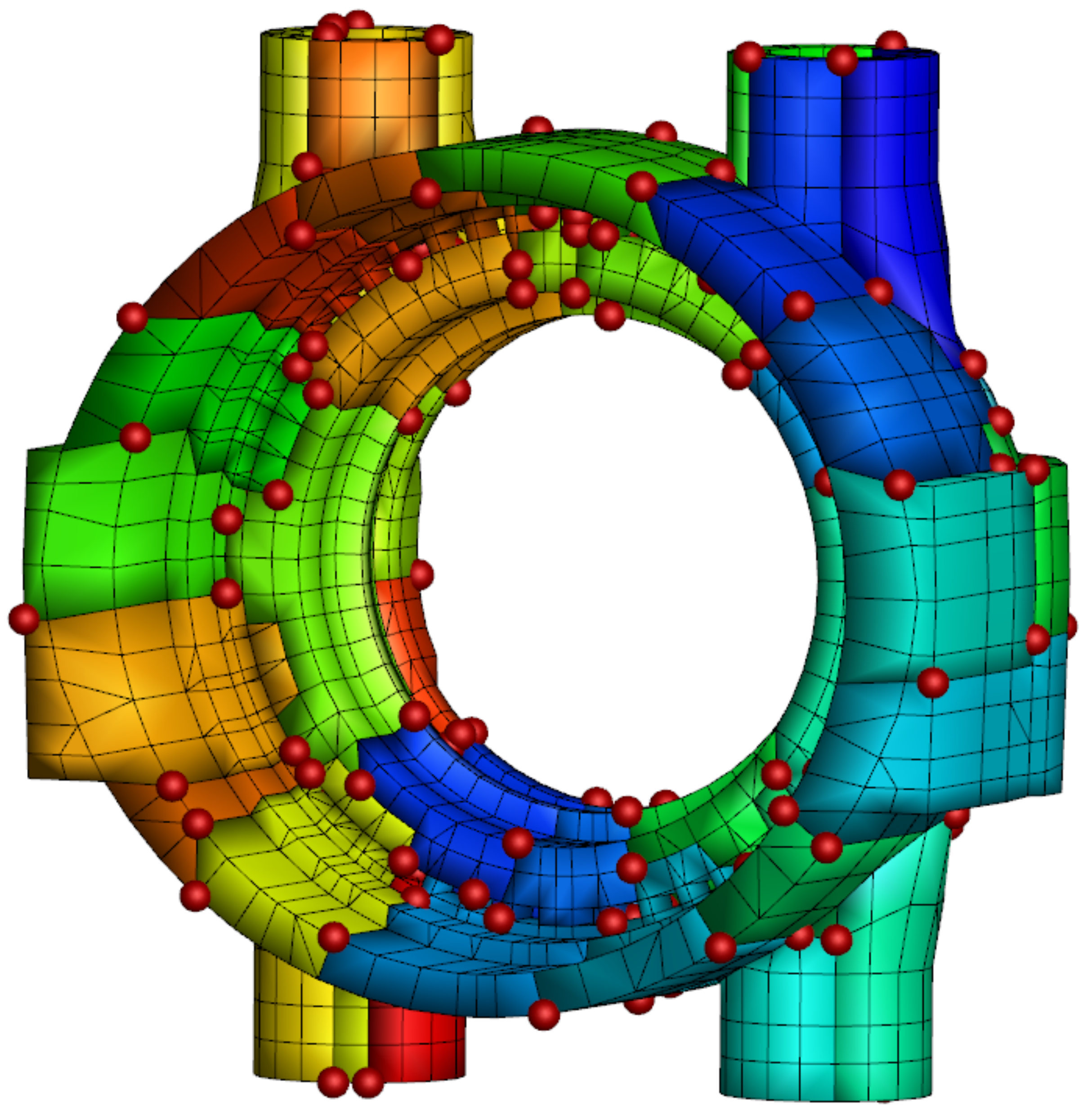}
\includegraphics[width=55mm]{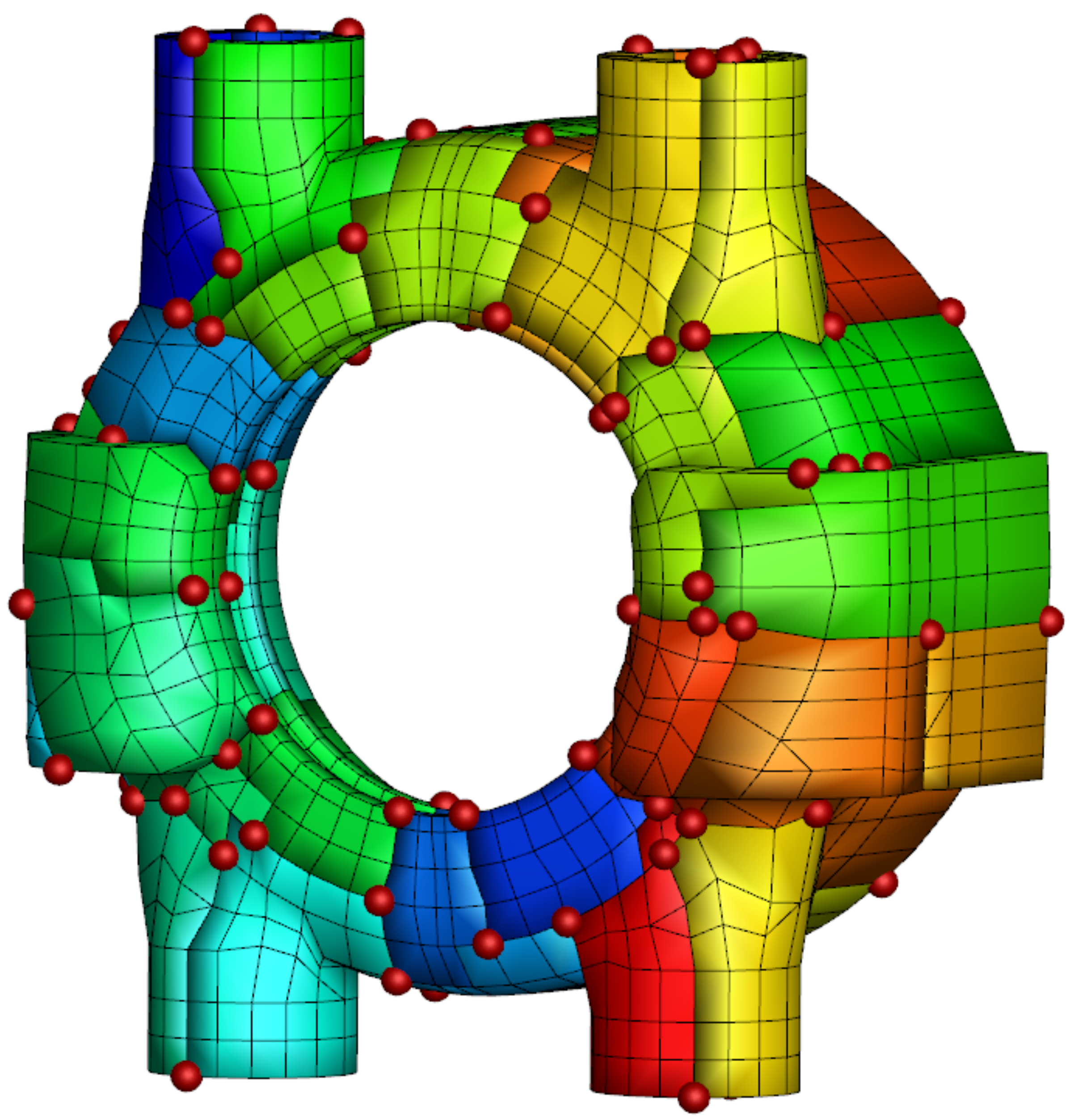}
\end{center}
\caption{Turbine nozzle problem, 36 subdomains, initial set of 218 corners selected by the full version of Algorithm~\ref{alg:selection} marked by balls.}
\label{fig:dyza_geo}
\end{figure}

\begin{figure}[htbp]
\begin{center}
\includegraphics[width=55mm]{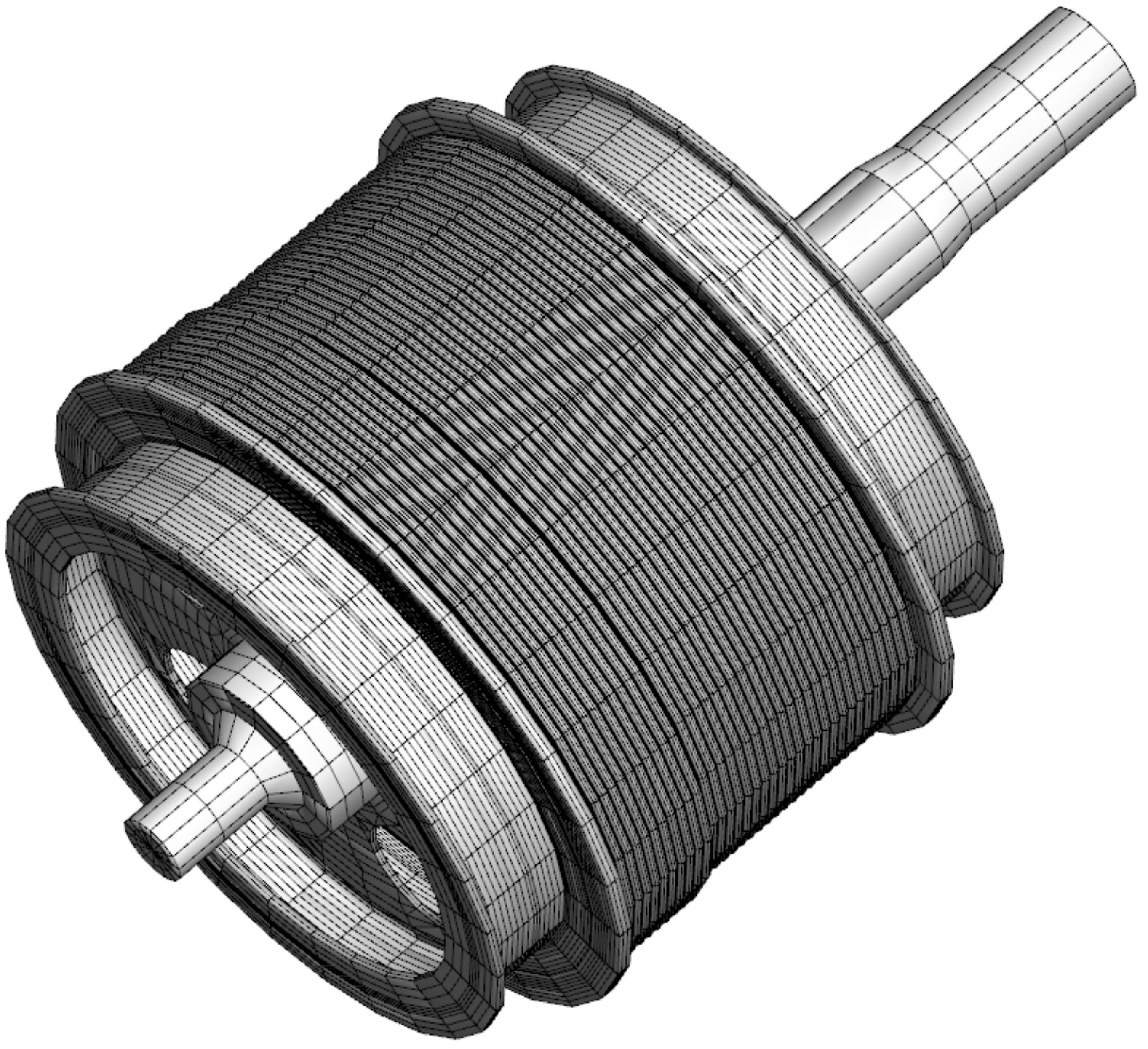}
\hskip 10mm
\includegraphics[width=55mm]{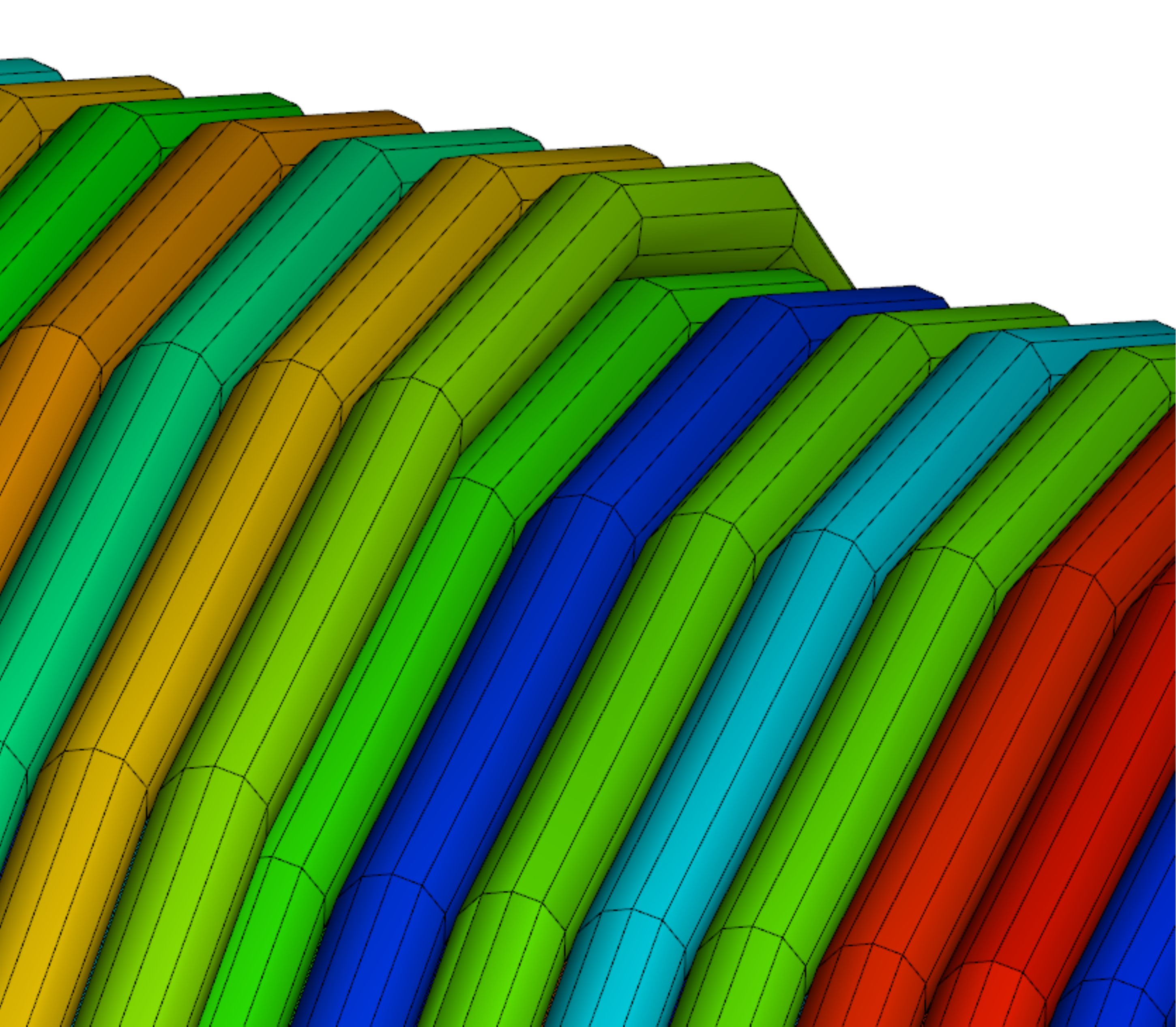}
\end{center}
\caption{Mine reel problem, finite element mesh (left) and a detail of the steel rope with division into subdomains (right).}
\label{fig:buben_geo}
\end{figure}

Three different industrial problems have been tested. 
The first one is a~problem of elasticity analysis of a turbine nozzle,
through which the steam enters the turbine blades (Figure~\ref{fig:dyza_geo}). The geometry is discretized using 2\,696 quadratic elements, which leads to 40\,254 unknowns. 
The second one is a~problem of elasticity analysis of a hip joint replacement which is loaded by pressure from human body weight. This mesh consists of 33\,186 quadratic elements resulting in 544\,734 unknowns.
Both meshes are divided into 36 subdomains by METIS. 
The turbine nozzle problem was computed using 12 processors, for hip joint replacement 36 processors were used.
The third problem is stress analysis of a mine reel loaded by its own weight and the weight of the steel wire rope (Figure~\ref{fig:buben_geo}). 
The mesh consists of 140\,816 quadratic elements and 1\,739\,211 unknowns.
It was divided into 1\,024 subdomains by METIS. Problem was computed using 32 processors.
Decomposition characteristics of the three industrial problems are summarized in Table \ref{tab:interf}.

Three algorithms for selecting the basic set of corners are tested: Algorithm~1 from Section \ref{sec:selection} referred to as \emph{full}, modified Algorithm~1 described in Remark~\ref{remark:favouring} in Section \ref{sec:selection} referred to as \emph{minimal}, and the edge-based algorithm mentioned in Section \ref{sec:selection}, inspired by \cite{Klawonn-2003-SCD} and described in \cite{Sistek-2008-FEM}, referred to as \emph{edge}. 
The number of PCG iterations was chosen as a measure of quality of the BDDC preconditioning. 
Numbers of the basic sets of corners obtained by the three algorithms for the three problems are recorded in Table \ref{tab:basic} and corresponding number of PCG iterations are summarized in Table \ref{tab:pcg}. For the two smaller problems (turbine nozzle and hip joint replacement), either constraints on corners only (referred to as \emph{C}), or constraints on corners and all averages (over all edges and faces) referred to as \emph{C+E+F} are tested. 
For the larger problem of mine reel, 
corner constraints alone turned out to be too weak to achieve a reasonable convergence and the results are marked as `n/a'. 
The edge-based algorithm did not work properly for hip joint replacement problem in the case of the basic set of corners only, 
so the results are missing too.

As the basic sets of corners selected by different algorithms have different numbers of corners, 
for a fair comparison of the algorithms we added more corners selected randomly from the interface to the smaller sets in order 
to achieve the same number of corners. 
Comparison of the algorithms using the same number of corner constraints is summarized in Table~\ref{tab:pcg_rand}. 

\begin{figure}[htbp]
\begin{center}
\includegraphics[width=60mm]{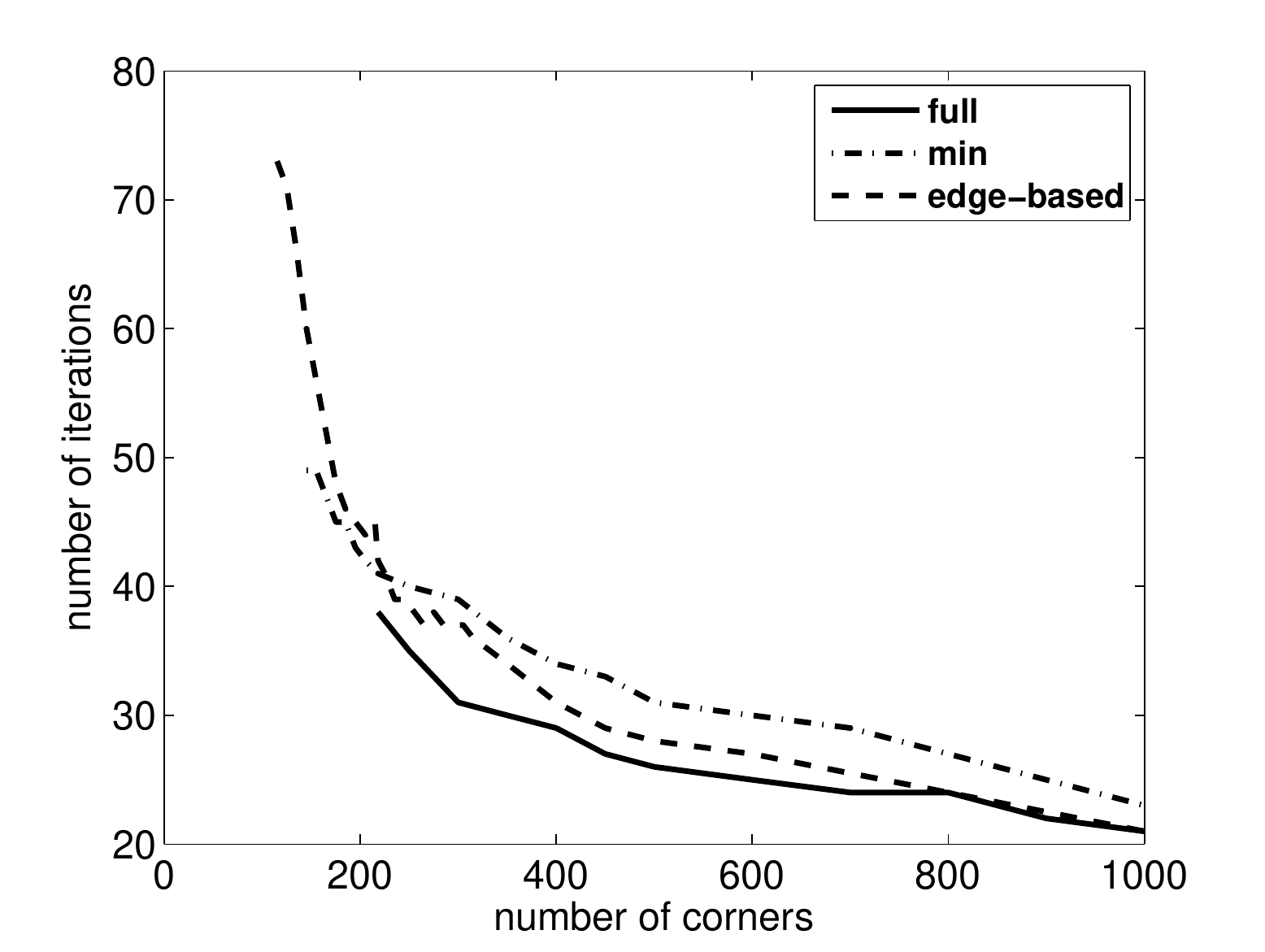}
\includegraphics[width=60mm]{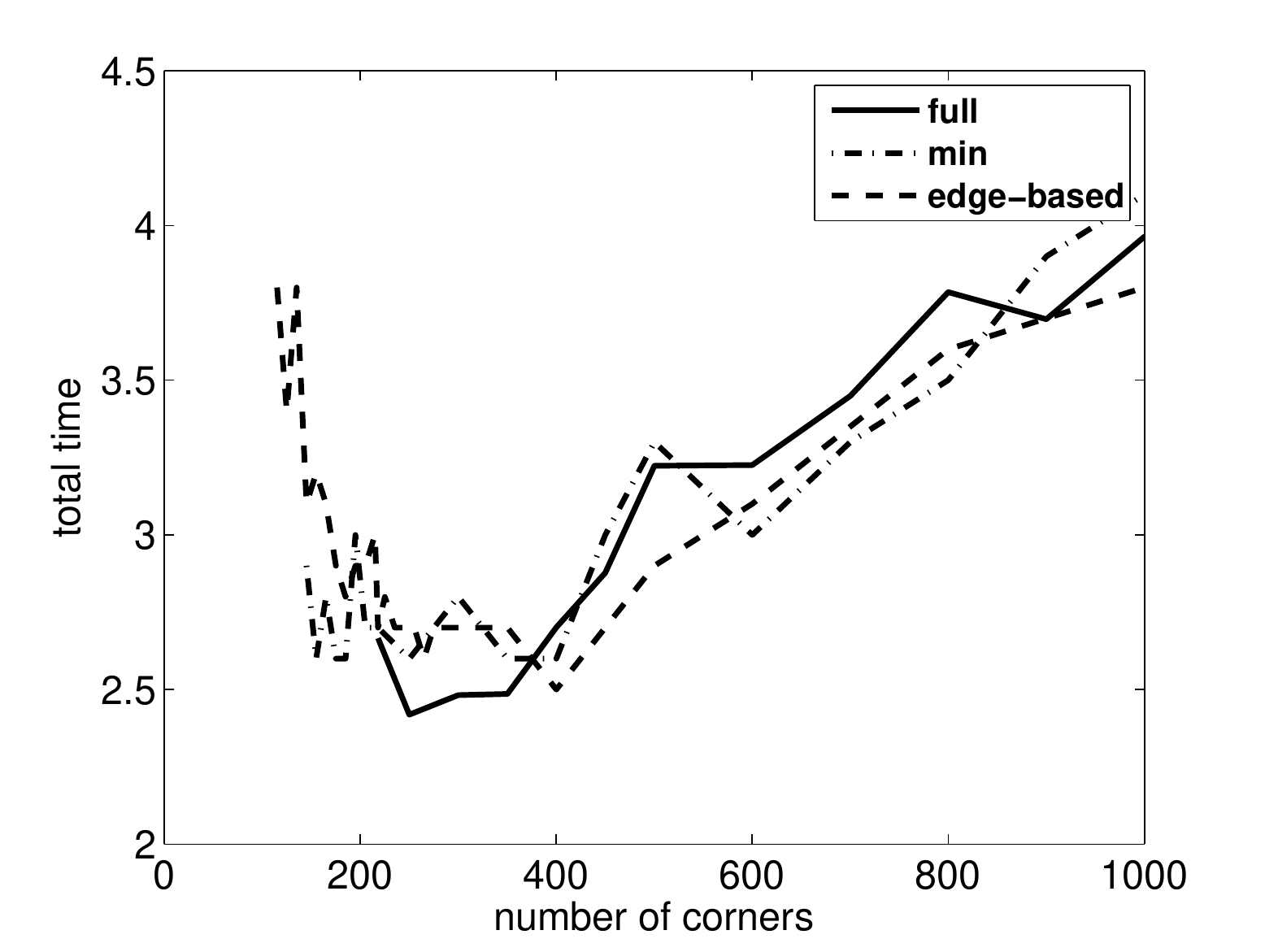}
\end{center}
\caption{Turbine nozzle problem, 36 subdomains, corner constraints only. 
Dependence of the number of iterations (left) and the total computational time (right) on the number of corner constraints. 
Full line - full version of the Algorithm~\ref{alg:selection}, dash-dotted line - minimalistic version, dashed line - the edge based algorithm.}
\label{fig:dyza}
\end{figure}

\begin{figure}[htbp]
\begin{center}
\includegraphics[width=60mm]{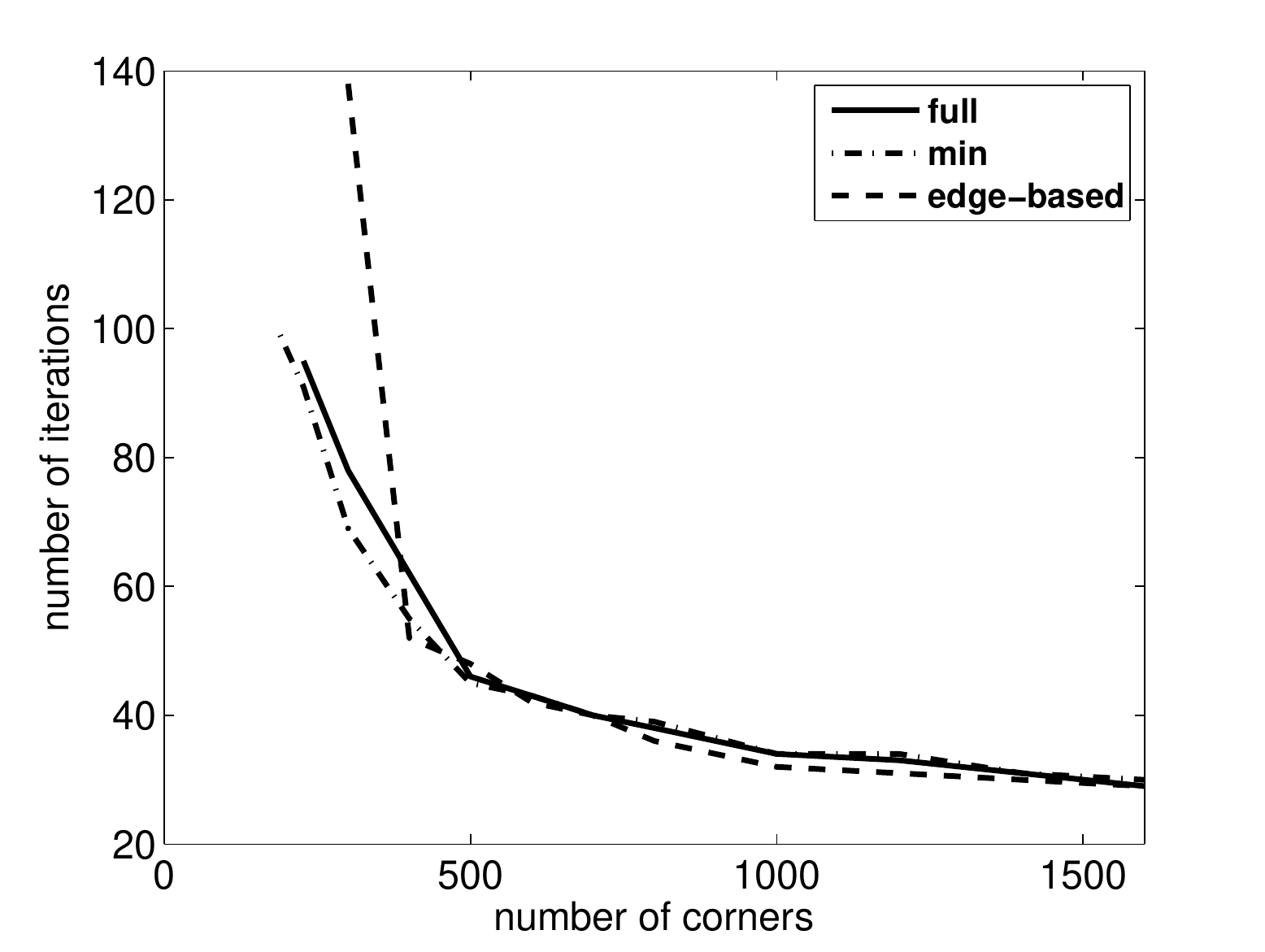}
\includegraphics[width=60mm]{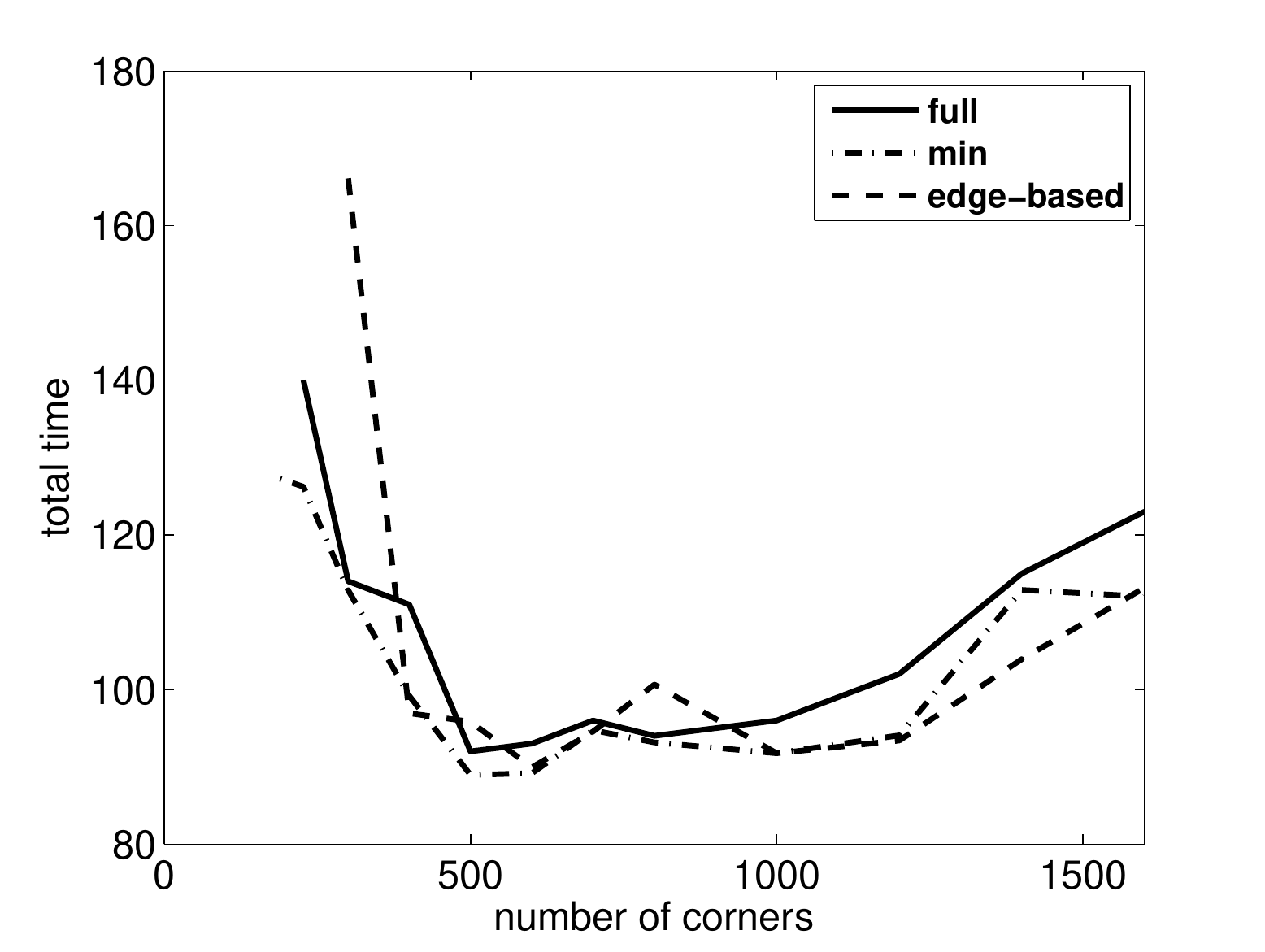}
\end{center}
\caption{Hip joint replacement problem, 36 subdomains, corner constraints only. 
Dependence of the number of iterations (left) and the total computational time (right) on the number of corner constraints. 
Full line - full version of the Algorithm~\ref{alg:selection}, dash-dotted line - minimalistic version, dashed line - the edge based algorithm.}
\label{fig:celekt}
\end{figure}

\begin{figure}[htbp]
\begin{center}
\includegraphics[width=60mm]{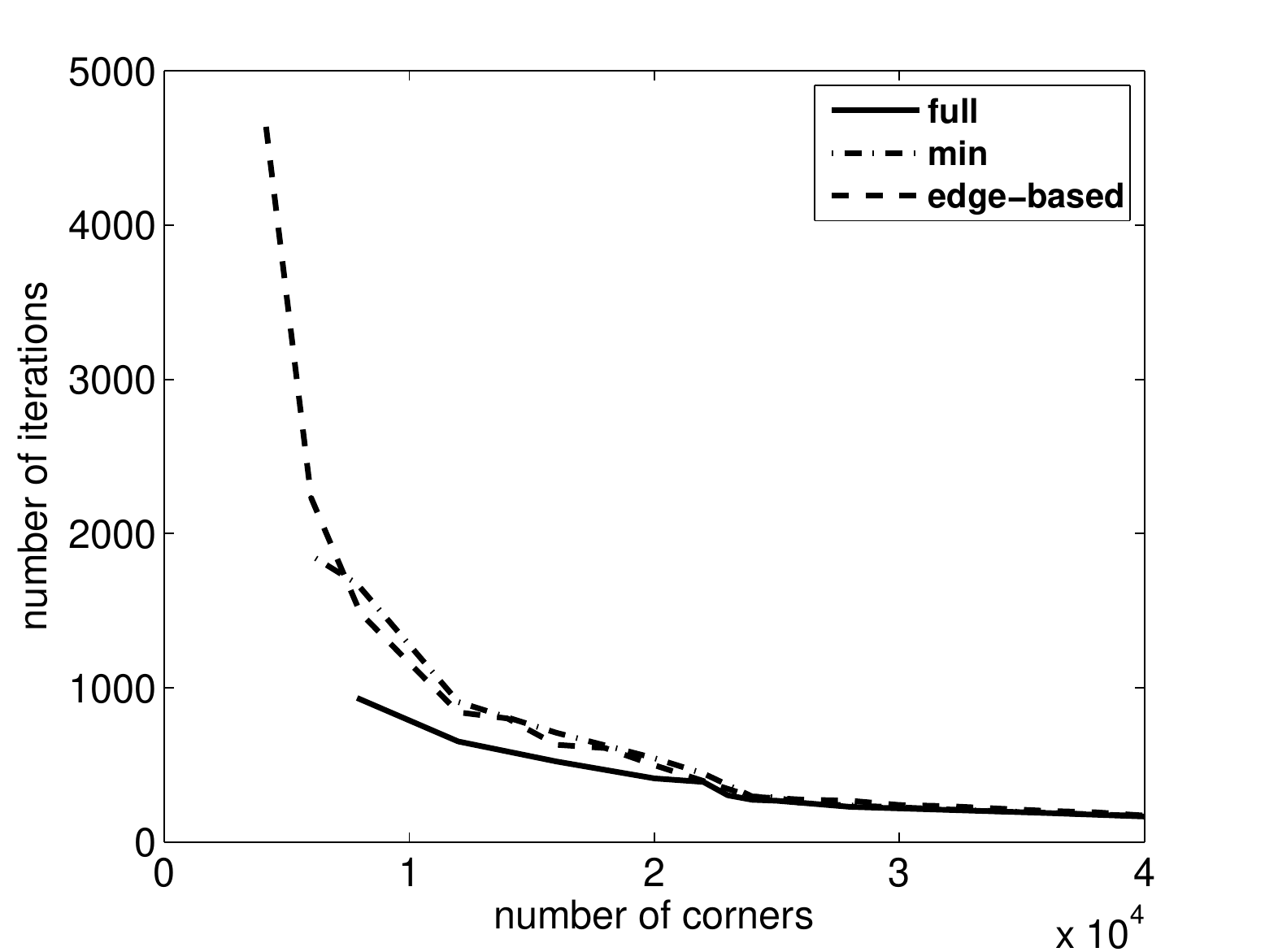}
\includegraphics[width=60mm]{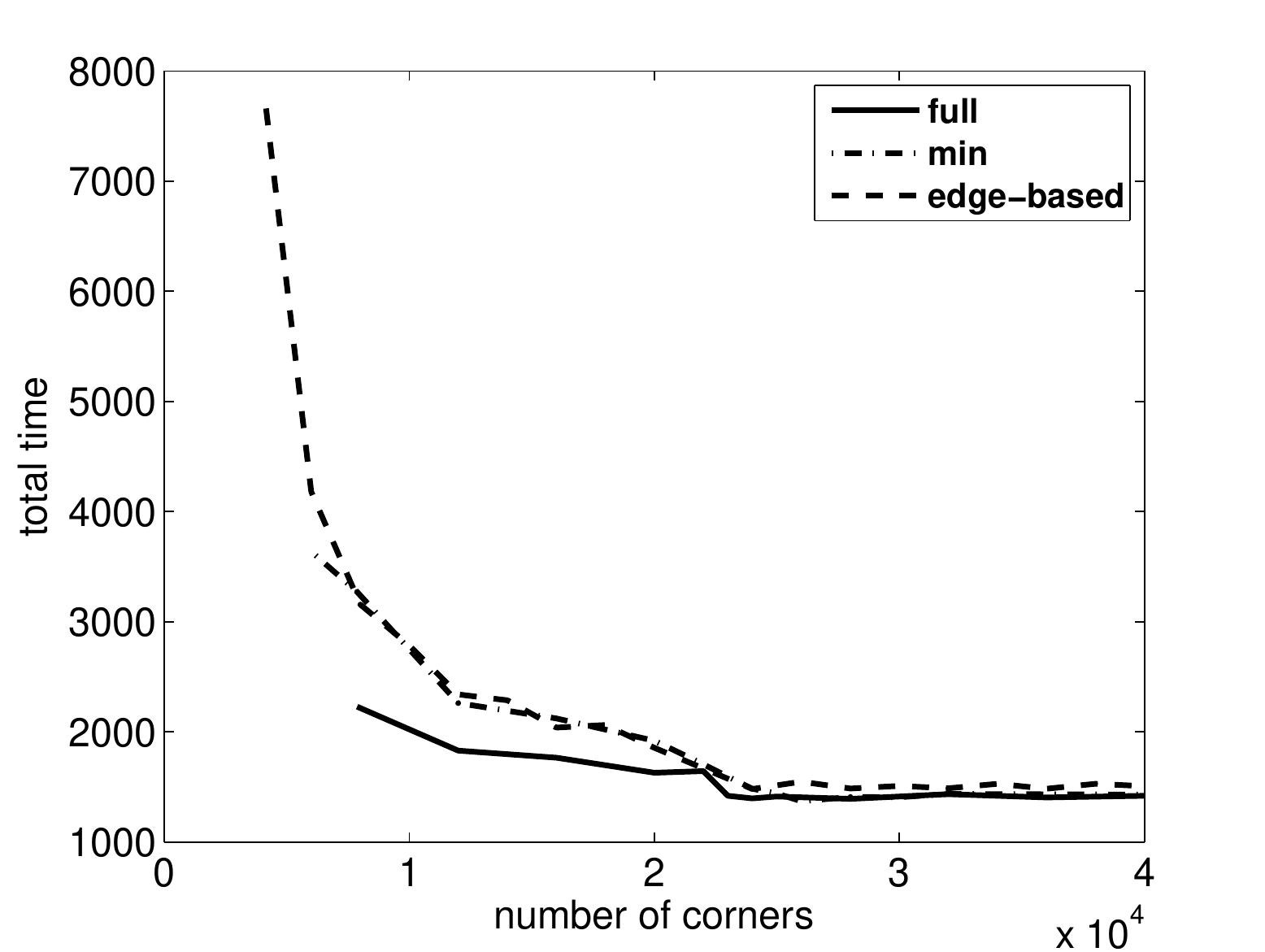}
\end{center}
\caption{Mine reel problem, 1024 subdomains, corner and all edge and face constraints. 
A~dependence of the number of iterations (left) and the total computational time (right) on the number of corner constraints. Full line - full version of the Algorithm~\ref{alg:selection}, 
dash-dotted line - minimalistic version, dashed line - the edge based algorithm.}
\label{fig:buben}
\end{figure}

Interesting results are obtained by adding more randomly selected interface nodes as corners to the basic set in order to improve convergence (see Figures \ref{fig:dyza} -- \ref{fig:buben} left): it seems that the initial choice of the basic set  influences the convergence properties even when many more randomly selected corners are added. Graphs on the right side of these figures show that the best computational time is achieved for higher numbers of corners than the basic sets for all problems tested and all algorithms for selecting the basic set used.

It can be observed especially on the most difficult problem of mine reel (Fig.~\ref{fig:buben}),
that the basic set of corners provided by the new algorithm in its full version is much more efficient than the basic sets provided by the earlier approaches and
considerably reduces the computational time.

\section{Conclusion}

It has been observed on a number of practical computations by the BDDC method, 
that the effort to find the minimal set of corners might be misleading and
selecting more corners often considerably improves the performance of the preconditioner and reduces the computational time.
This behaviour can be explained for problems with complex interface by position of selected corners,
which may be optimal with respect to one pair of subdomains, but may lead to poorly conditioned problems for other subdomains. 
As a consequence, by favouring already selected corners, 
these subdomains are not given the freedom to select corners optimally distributed for their own fixation.

This has been the main motivation for presenting a new approach to selecting the basic set of corners,
which is proposed in Section~\ref{sec:selection}.
It attempts to combine advantages of previous algorithms, and it is based on selection of corners independently for each face, 
so it can be naturally parallelized.
It does not aspire to minimize the number of selected corners that assure the invertibility of all problems in BDDC
and typically produces a larger initial set of coarse nodes than the other algorithms. 
We have seen this to be beneficial for all performed computations.

Numerical experiments on three industrial problems show that for basic sets of corners,
this approach gives better results than the other two algorithms used for comparison in all three tested problems. 
When more corners are added, 
better results are obtained in two of the problems (turbine nozzle and mine reel) and comparable results in the third case (hip joint replacement).

We are aware that for very large problems the solution of the coarse problem might eventually dominate the computation and another approach
than a (parallel) direct solver could be necessary.
In such cases, multilevel extension of the BDDC method (e.g. \cite{Mandel-2008-MMB}) seems to be a promising way.
However, we observed even for the largest test problem of the mine reel,
that we did not reach this computational bottleneck when adding more corners into the coarse problem,
and the curve of computational time with respect to the number of corners was still decreasing.
The expected bottleneck is also pushed farther by the everlasting advances in parallel direct solvers.

\section*{Acknowledgement}

We are grateful to Jan Le\v{s}tina (Vamet Ltd.) for providing the problems of turbine nozzle and mine reel for testing.
We would also like to thank to Jan Mandel and Bed\v rich Soused\'\i k for fruitful discussions about this topic and BDDC method in general.
Jakub \v{S}\'{\i}stek would like to acknowledge the help of Lan Vu with improvement of the algorithm for disconnected subdomains.
This research has been supported by the Czech Science Foundation under grant
GA CR 106/08/0403, by the Grant Agency of the Academy of Sciences of the CR under grant IAA100760702, and by
projects MSM 6840770001 and MSM 6840770010. It has also been supported
by Institutional Research Plans AV0Z 10190503 and AV0Z 20760514.

\bibliographystyle{elsarticle-harv}
\bibliography{/home/sistek/denver/bddc/bibliography/bddc}

\newpage

\begin{table}[htp]
\begin{center}
\begin{tabular}
[c]{|c||c|r|r|r|r|r|}\hline
problem  & subs.  & vertices & edges & faces & intf. nodes & all nodes\\\hline\hline
Turbine nozzle & 36 & 6  & 60  & 101  & 2 714 & 13 418  \\\hline
Hip replacement & 36 & 1  & 19  & 78   & 9 222 & 181 578  \\\hline
Mine reel      & 1 024 & 2 451  & 1 209  & 4 164  & 117 113 & 579 737  \\\hline
\end{tabular}
\end{center}
\caption{Decomposition characteristics of the tested problems.}
\label{tab:interf}
\end{table}

\begin{table}[htp]
\begin{center}
\begin{tabular}
[c]{|c||r|r|r|}\hline
problem  & full  & min & edge \\\hline\hline
Turbine nozzle & 218 & 145  & 115 \\\hline
Hip replacement & 227 & 189  & 66 \\\hline
Mine reel      & 7 864 & 6 183  & 4 152  \\\hline
\end{tabular}
\end{center}
\caption{Number of corners in the basic set selected by different algorithms.}
\label{tab:basic}
\end{table}

\begin{table}[htp]
\begin{center}
\begin{tabular}
[c]{|c||r|r|r|r|r|r|}\hline 
 & \multicolumn{3}{|c|}{C} & \multicolumn{3}{|c|}{C+E+F} \\
  & full  & min & edge & full  & min &  edge \\\hline\hline

Turbine nozzle & 38 & 49  & 73  & 24  & 27 & 29  \\\hline
Hip replacement & 95 & 99  & n/a  & 50   & 52 & n/a \\\hline
Mine reel      & n/a  & n/a   & n/a   & 935  & 1 841 & 4 637  \\\hline
\end{tabular}
\end{center}
\caption{Number of PCG iterations needed for convergence for different algorithms of selecting the basic set of corners and different constraint type.}
\label{tab:pcg}
\end{table}

\begin{table}[htp]
\begin{center}
\begin{tabular}
[c]{|c||r|r|r|r|r|r|}\hline 
 & \multicolumn{3}{|c|}{C} & \multicolumn{3}{|c|}{C+E+F} \\
  & full  & min & edge & full  & min &  edge \\\hline\hline

Turbine nozzle & 38 & 41  & 42  & 24  & 25 & 26  \\\hline
Hip replacement & 95 & 91  & $>$ 138  & 50   & 50 & 61  \\\hline
Mine reel      & n/a  & n/a   & n/a   & 935  & 1 674 & $\approx$1 800  \\\hline
\end{tabular}
\end{center}
\caption{Number of PCG iterations needed for convergence for different algorithms of selecting the basic set of corners and different constraint type. For every problem, different basic sets were completed to the same number of corners by adding randomly selected corners.}
\label{tab:pcg_rand}
\end{table}

\end{document}